\documentclass[epic,11pt]{article}
\usepackage{amsfonts,amsmath,amscd,amssymb,amsthm}
\usepackage{fullpage}
\usepackage[pdftex]{pict2e}

\usepackage[alphabetic,initials]{amsrefs}
\usepackage[all]{xy}
\usepackage{hyperref} 
\usepackage{mathrsfs}
\input xy%
\usepackage{fullpage}
\usepackage[pdftex]{pict2e}
\usepackage[dvipsnames]{xcolor}
\def\0{{\bar 0}}
\def\1{{\bar 1}}
\def\mo{{\operatorname{Mor}}}

\def\Z{{\mathbb Z}}
\def\X{{\mathbb X}}
\def\F{{\mathbb F}}

\def\B{{\mathbb B}}

\def\ua{{\operatorname{a  }}}
\def\ur{{\operatorname{r  }}}
\def\ud{{\operatorname{d}}}

\def\red{\operatorname{red}}

\def\span{\operatorname{span\;}}

\def\str{{\operatorn<ame{str}}}
 
\def\trans{{\operatorname{transpose}}} 

\def\dist{{\operatorname{dist}}}

\def\iso{{\operatorname{iso}}}

\newcommand{\ttk}{\mathtt{k}}

\newcommand{\lst}{{\stackrel{\leftarrow}{\gl}}}

\newcommand{\ghs}{{\gth(\gs)}}

\newcommand{\bA}{{\bf a}}

\newcommand{\diag}{diag}

\newcommand{\Lgo}{L(\stackrel{{\rm _o}}{\fg})}
\newcommand{\Lgh}{\widehat{L}(\stackrel{{\rm _o}}{\fg})}

\newcommand{\itema}{\item[{{\rm$($a$)$}}]}
\newcommand{\itemb}{\item[{{\rm$($b$)$}}]}
\newcommand{\itemc}{\item[{{\rm$($c$)$}}]}
\newcommand{\itemd}{\item[{{\rm$($d$)$}}]}
\newcommand{\iteme}{\item[{{\rm$($e$)$}}]}

\newcommand{\da}{\downarrow}
\newcommand{\pa}{\uparrow}
\newcommand{\itemo}{\item[{}]}
\newcommand{\noi}{\noindent}

\newcommand{\ga}{\alpha}
\newcommand{\gb}{\beta}
\newcommand{\gc}{\gamma}

\newcommand{\Gl}{\Lambda}
\newcommand{\Gd}{\Delta}

\newcommand{\gd}{\delta}

\newcommand{\gs}{\sigma}

\newcommand{\gt}{\tau}
\newcommand{\gz}{\zeta}

\newcommand{\gl}{\lambda}

\newcommand{\gr}{\rho}
\newcommand{\gk}{\kappa}
\newcommand{\gep}{\epsilon}
\newcommand{\gth}{\theta}
\newcommand{\op}{\oplus}
\newcommand{\A}{\mathbb A}

\def\str{{\operatorname{str}}}

\def\Im{{\operatorname{Im}\;}}

\newcommand{\fg}{\mathfrak{g}}\newcommand{\fgl}{\mathfrak{gl}}
\newcommand{\fsl}{\mathfrak{sl}}\newcommand{\osp}{\mathfrak{osp}}

\newcommand{\fr}{\mathfrak{r}}

\newcommand{\fh}{\mathfrak{h}}

\newcommand{\fb}{\mathfrak{b}}

\newcommand{\fa}{\mathfrak{a}}

\newcommand{\fk}{\mathfrak{k}}

\newcommand{\fB}{\mathfrak{B}}
\newcommand{\fG}{\mathfrak{G}}
\newcommand{\fS}{\mathfrak S}
\newcommand{\fT}{\mathfrak T}

\newfont{\eufm}{eufm10 scaled\magstep1}

 \newcommand{\ti}{\times}

\newcommand{\bcu}{\bigcup}
\newcommand{\bop}{\bigoplus}
\newcommand{\bsk}{\backslash}

\newcommand{\cO}{\mathcal{O}}

\newcommand{\cB}{\mathcal{B}}

\newcommand{\cS}{\mathcal{S}}

\newcommand{\cH}{\mathcal{H}}

\newcommand{\cW}{\mathfrak{W}}

\newcommand{\Sh}{{\bf{Sh}}}
\newcommand{\ey}{\end{eqnarray}}
\newcommand{\by}{\begin{eqnarray}}
\newcommand{\nn}{\nonumber}

\newcommand{\bco}{\begin{conjecture}}
\newcommand{\ba}{\begin{alg}}
\newcommand{\ea}{\end{alg}}
\newcommand{\eco}{\end{conjecture}}
\newcommand{\bpf}{\begin{proof}}
\newcommand{\epf}{\end{proof}}
\newcommand{\bt}{\begin{theorem}}

\newcommand{\et}{\end{theorem}}

\newcommand{\br}{\begin{rem}}
\newcommand{\er}{\end{rem}}
\newcommand{\brs}{\begin{rems}}
\newcommand{\ers}{\end{rems}}
\newcommand{\bi}{\begin{itemize}}
\newcommand{\ei}{\end{itemize}}
\newcommand{\bl}{\begin{lemma}}
\newcommand{\bsul}{\begin{sublemma}}
\newcommand{\esul}{\end{sublemma}}
\newcommand{\bp}{\begin{proposition}}
\newcommand{\be}{\begin{equation}}
\newcommand{\bc}{\begin{corollary}}
\newcommand{\bexs}{\begin{examples}}
\newcommand{\eexs}{\end{examples}}
\newcommand{\bexa}{\begin{example}}
\newcommand{\eexa}{\end{example}}
\newcommand{\bex}{\begin{exercise}}
\newcommand{\eex}{\end{exercise}}
\newcommand{\btab}{\begin{tab}}
\newcommand{\etab}{\end{tab}}
\newcommand{\el}{\end{lemma}}
\newcommand{\ep}{\end{proposition}}
\newcommand{\ee}{\end{equation}}
\newcommand{\ec}{\end{corollary}}
\newcommand{\Bc}{\begin{center}}
\newcommand{\Ec}{\end{center}}
\newcommand{\bh}{\begin{hyp}}
\newcommand{\eh}{\end{hyp}}
\newcommand{\bhs}{\begin{hyps}}
\newcommand{\ehs}{\end{hyps}}
\newcommand{\bd}{\begin{dfn}}
\newcommand{\ed}{\end{dfn}}

\newcommand{\bn}{\begin{notn}}
\newcommand{\en}{\end{notn}}
\newcommand{\ogd}{\overline{\delta}}

\newcommand{\fgm}{\stackrel{{\rm o}}{\fg}^{-\gth}}

\newcommand{\fgo}{\stackrel{{\rm o}}{\fg}}
\newcommand{\Fbo}{\stackrel{{\rm o}}{\cB}}
\newcommand{\pigo}{\stackrel{{\rm o}}{\Pi}}
\newcommand{\fbo}{\stackrel{\rm o}{\fb}}{}
\newcommand{\pigov}{\stackrel{\rm o}{\Pi^\vee}}

\newcommand{\fgt}{\stackrel{{\rm o}}{\fg}^\gth}

\newcommand{\fho}{\stackrel{\rm o}{\fh}}
\newcommand{\fhs}{\stackrel{{\rm _o}}{\fh}^*}

\usepackage[aligntableaux=center]{ytableau}
\usepackage{tabularray}

\begin{document}
\title{Table of Contents}


\newtheorem*{bend}{Dangerous Bend}

\newtheorem{thm}{Theorem}[section]
\newtheorem{hyp}[thm]{Hypothesis}
 \newtheorem{hyps}[thm]{Hypotheses}
\newtheorem{notn}[thm]{Notation}

  \newtheorem{rems}[thm]{Remarks}

\newtheorem{conjecture}[thm]{Conjecture}
\newtheorem{theorem}[thm]{Theorem}
\newtheorem{theorem a}[thm]{Theorem A}
\newtheorem{example}[thm]{Example}
\newtheorem{examples}[thm]{Examples}
\newtheorem{corollary}[thm]{Corollary}
\newtheorem{rem}[thm]{Remark}
\newtheorem{lemma}[thm]{Lemma}
\newtheorem{sublemma}[thm]{Sublemma}
\newtheorem{cor}[thm]{Corollary}
\newtheorem{proposition}[thm]{Proposition}
\newtheorem{exs}[thm]{Examples}
\newtheorem{ex}[thm]{Example}
\newtheorem{exercise}[thm]{Exercise}
\numberwithin{equation}{section}%
\setcounter{part}{0}
\newcommand{\drar}{\rightarrow}
\newcommand{\lra}{\longrightarrow}
\newcommand{\rra}{\longleftarrow}
\newcommand{\dra}{\Rightarrow}
\newcommand{\dla}{\Leftarrow}
\newcommand{\rl}{\longleftrightarrow}

\newtheorem{Thm}{Main Theorem}


\newtheorem*{thm*}{Theorem}
\newtheorem{lem}[thm]{Lemma}
\newtheorem*{lem*}{Lemma}
\newtheorem*{prop*}{Proposition}
\newtheorem*{cor*}{Corollary}
\newtheorem{dfn}[thm]{Definition}
\newtheorem*{defn*}{Definition}
\newtheorem{notadefn}[thm]{Notation and Definition}
\newtheorem*{notadefn*}{Notation and Definition}
\newtheorem{nota}[thm]{Notation}
\newtheorem*{nota*}{Notation}
\newtheorem{note}[thm]{Remark}
\newtheorem*{note*}{Remark}
\newtheorem*{notes*}{Remarks}
\newtheorem{hypo}[thm]{Hypothesis}
\newtheorem*{ex*}{Example}
\newtheorem{prob}[thm]{Problems}
\newtheorem{conj}[thm]{Conjecture}

\title{The Weyl groupoid in Type A, Young diagrams and Borel subalgebras}
\author{Ian M. Musson
 \\Department of Mathematical Sciences\\
University of Wisconsin-Milwaukee\\ email: {\tt
musson@uwm.edu}}
\maketitle
\section*{ \color{white} phantom}
   \textit{
A little cylinder of paper fell from the top of my locker. As I unrolled it I saw  black calligraphy crawling across it like a spider.}\\
   \vspace{.2cm}
   - Alexandra Adornetto \cite{A}
\begin{abstract} Let $\ttk$ be an algebraically closed field of characteristic zero.  Let ${\stackrel{{\rm o}}{\fg}}$ be the Lie superalgebra $\fsl(n|m)$ and let  $\mathfrak{W}$  be the 
Weyl groupoid 
introduced by 
 Sergeev and Veselov \cite{SV2} using the root system of ${\stackrel{{\rm o}}{\fg}}$.
An important subgroupoid 
$\mathfrak T_{iso}$ of $\cW$ has base 
$\Gd_{iso}$, the set of all the isotropic roots.  
  In \cite{SV101}, motivated by deformed quantum Calogero-Moser problems \cite{SV1},  the same authors considered an action of $\mathfrak{W}$ on $\ttk^{n|m}$   
depending on a parameter $\gk$.  
In the case $m>n$, with $m,n $ relatively prime  and $\gk=-n/m$ we study a particular infinite orbit 
of $\mathfrak T_{iso}$ 
with some special properties. This orbit, thought of as a directed graph is isomorphic to the graph of an orbit for the action of $\mathfrak T_{iso}$ 
on certain Borel subalgebras of  the affinization $\Lgh$ of  $\fgo$. 

The underlying reason for this graph isomorphism is that both have combinatorics which can be described using  Young diagrams and tableaux 
drawn on the surface of a rotating cylinder with circumference $n$ and length $m$. 
Allowing the  cylinder to rotate produces an infinite orbit. 
This leads to a third graph which is isomorphic to the other two.
\end{abstract}


\section{Introduction}\label{itr} 
Throughout  $[k]$ denotes  the set of the first $k$ positive integers and iff means if and only if. 
We show that the orbits arising from certain groupoid actions are isomorphic as directed graphs. 
For a set  $X$, the symmetric groupoid  $\mathfrak \fS(X)$ on $X$ has base all subsets of $X$ and  for $U,V\subseteq X$, the set of morphisms $\mo(U,V)$ is the set of all bijections from $U$ to $V$.  A groupoid $\fG$ {\it acts} on $X$ if there is 
a functor  
\be \label{cwr}F:\fG\lra \mathfrak S(X).\ee
 The groupoid actions  arise in essentially three different ways.  First, given positive integers $m, n$ with $m\ge n$, let $X$ be the set of Young diagrams that fit inside a $m\ti n$ rectangle {\bf R}.  Adding/deleting boxes to/from a diagram are operations on $X$ which are not always defined, and this observation leads to the construction of a groupoid $\mathfrak T_{iso}(X)$ which acts on $X$ in a natural way. 
Let ${\stackrel{{\rm o}}{\cB}}$ be the set of Borel subalgebras $\fbo$ of $\fgl(n|m)$ such that 
$\fbo_0$ is the Borel subalgebra of $\fg$ consisting of two diagonal blocks of 
upper triangular matrices.  There is a natural bijection $X\rl {\stackrel{{\rm o}}{\cB}}$ and this induces an action of $\mathfrak T_{iso}$   on ${\stackrel{{\rm o}}{\cB}}$ by odd reflections.  
 For details including some notation that we postpone, see Subsection \ref{yd}. 
Now informally, the action  of $\mathfrak T_{iso}$ on $X$ can be continued indefinitely by deleting the first row of a Young diagram $\gl $ whenever $\gl_1=m$ and keeping track of the number of times a row is deleted.
In more detail, if $\gl =(\gl_1, \gl_2, \ldots \gl_n)\subseteq$ {\bf R}, 
and
$\gl_1 = m$, 
let   $\overline \gl =  (\gl_2, \ldots,\gl_n,0 )$  be the partition obtained from removing the first part from $\gl$.  
Let $\sim$ be the smallest equivalence relation on 
$X\ti \Z$ such that $(\gl,k-1)\sim (\overline \gl,k)$ if $\gl_1=m$.   There is a compatible 
equivalence relation on 
${\stackrel{{\rm o}}{\cB}}\ti \Z.$  The action of $\mathfrak T_{iso}$  can be extended to the set of equivalence classes  $[X \ti \Z]$ and   $[{\stackrel{{\rm o}}{\cB}} \ti \Z]$. 
We show in Corollary \ref{sz} that there are functors $F, B$ that fit into a commutative diagram. 
\be\label{sez}
\xymatrix@C=2pc@R=1pc{
\mathfrak  T_{iso}\ar@{<->}_= [dd]&&
\fS[X \ti \Z] \ar@{<-}_F[ll]: \ar@{<->}^{\cong}[dd]&\\ \\
\mathfrak  T_{iso} \ar@{-}[rr]_B &&
\fS[{\stackrel{{\rm o}}{\cB}} \ti \Z]\ar@{<-}[ll] &}
\ee
Suppose we have a functor $G:\mathfrak  T_{iso}\lra \mathfrak S(V)$ and $v, w\in V$. Write  
\be \label{ags} 
v \stackrel{{\ga}}{\longrightarrow } w\ee
 to mean that 
$v\in G(\ga)$ and $G(\gr_\ga)(v) =w.$ We have a graph  with vertex set $V$ and labelled arrows as in \eqref{ags}.  The connected component of this graph containing $v$ is the orbit of $v$.  We denote it by $\cO(G,v).$ 
Let ${\stackrel{{\rm o}}{\fb}}{}^\dist$ be the distinguished Borel subalgebra of $\fgo$ and denote the empty Young diagram $0^n$ by $\emptyset$.
From the commutative diagram \eqref{sez} we obtain an isomorphism of directed labelled graphs 
\be  \label{fx} \cO(F,[\emptyset, 0])  \lra \cO(B,[{\stackrel{{\rm o}}{\fb}}{}^\dist,0])\ee
such that the equivalence classes $[\emptyset, 0]$ and $[{\stackrel{{\rm o}}{\fb}}{}^\dist, 0]$ correspond. 
\\\\
\noi 
The Lie superalgebras $\fgl(n|m)$ and  $\fsl(n|m)$ have the same root system $\Gd$
and there is an obvious bijection between Borel subalgebras of  $\fgl(n|m)$  and $\fsl(n|m)$.  
It will be convenient to work with
$\fgo= \fgl(n|m)$ until Section \ref{akmls}
where $\fgo=\fsl(n|m)$ with $m>n$ or $\fgl(n|n)$.  The groupoid $\mathfrak T_{iso}(X)$ from the previous paragraph  is isomorphic to a subgroupoid $\mathfrak T_{iso}$ of the Weyl groupoid 
$\cW$ defined in \cite{SV2} using the root system of $\fgo$.  
In  \cite{SV101} Sergeev and Veselov defined 
an action of the   $\mathfrak{W}$  on $\ttk^{n|m}$.  We reformulate their definition so that 
 $\mathfrak{W}$  acts on $\Z^{n|m}$ via a functor
 $$SV:\mathfrak {W}\lra \mathfrak S(\Z^{n|m}).$$ The restriction of $SV$
to $\mathfrak  T_{iso}$ is also denoted $SV$.  We consider the $\mathfrak  T_{iso}$-orbit of 
$\Gl_0 \in \Z^{n|m}$ as defined in \eqref{e7}.
\\\\
The Lie superalgebra  $\fgo$ has only finitely many conjugacy classes of Borel subalgebras, but the affinization  $\Lgh$ has infinitely many Borels which arise in the following way.
The process of affinization produces subalgebras ${\stackrel{\rm o}{\fg}}(k)$ for $k\in\Z,$ all  isomorphic to $\fgo$.  In addition for each shuffle $\gs$ there is a Borel subalgebra ${\stackrel{{\rm o}}{\fb}}{}(\gs,k)$ 
of ${\stackrel{\rm o}{\fg}}(k)$ as in \cite{M101}  Chapter 3 and each subalgebra ${\stackrel{{\rm o}}{\fb}}{}(\gs,k)$ extends to a unique Borel subalgebra of $\Lgh$. 
Also by affinization any Dynkin-Kac diagram of ${\stackrel{\rm o}{\fg}}(k)$ is  embedded in a corresponding diagram for  $\Lgh$ by adding an 
extending node.  Under suitable  conditions, by deleting a different node we obtain 
the next subalgebra ${\stackrel{\rm o}{\fg}}(k+1)$ together with a Borel subalgebra 
${\stackrel{{\rm o}}{\fb}}{}(\overline\gs,k+1)$.  
The Borels ${\stackrel{{\rm o}}{\fb}}{}(\gs,k)$ and ${\stackrel{{\rm o}}{\fb}}{}(\overline\gs,k+1)$ have the same extension to $\Lgh$.  
Let $\cB$ be the set of Borel subalgebras of  $\Lgh$ that are extensions of  Borels in the subalgebras ${\stackrel{\rm o}{\fg}}(k)$ for  $k\in \Z.$  
There is a functor  $B:\mathfrak  T_{iso}\lra \fS(\cB)$. Our main result is the following.

\bt \label{iir} \bi \itemo\itema There is an isomorphism of directed labelled graphs 
\be  \label{fy} \cO(F,(\emptyset, 0))  \lra \cO(B,\fb)\ee
where $\fb$ is extended from the distinguished subalgebra of ${\stackrel{\rm o}{\fg}}(0).$
\itemb If $m, n$ are coprime then the map  $x$  from 
Theorem \ref{rjt} 
induces an isomorphism of directed labelled graphs 
$$x:\cO(F,(\emptyset, 0))  \lra \cO(SV,\Gl_0)$$
such that $x(\emptyset, 0)  \lra \Gl_0$. 
\ei
\et \noi   In Theorem \ref{srb}, we show there is a bijection  
$\cB \rl [{\stackrel{{\rm o}}{\cB}} \ti \Z]$  which is compatible with the groupoid actions.   Part (a) of the Theorem follows from this and 
\eqref{fx}.  We also need a certain rooted directed subgraph of $\cO(B,[{\stackrel{{\rm o}}{\fb}}{}^\dist,0])$.  This subgraph $\cO^+(B,[{\stackrel{{\rm o}}{\fb}}{}^\dist,0])$ is obtained from the root $[{\stackrel{{\rm o}}{\fb}}{}^\dist,0]$ by applying (compositions of) positive morphisms.  
\\ \\
The three graphs 
have their origins in diverse areas: Young diagrams, Borel subalgebras and orbits of the Weyl groupoid under the action defined \cite{SV101}.  Hopefully this work serves to elucidate the similarites and differences in the 
combinatorics of these areas.
This paper is organized as follows. Section \ref{fc} contains preliminary results on Young diagrams, groupoids and Borel subalgebras. This Section includes a ``baby version",  \ref{bby} of diagram \ref{sez}.
\\ \\
In Section \ref{APG}   we explain how the action of $\mathfrak  T_{iso}$ can be extended to an action on a set of equivalence classes. 
Orbits for the action of $\mathfrak  T_{iso}$  on $\Z^{n|m}$ 
are studied in 
Section \ref{affc}.   
In Section \ref{akmls},  we discuss 
Borel subalgebras of 
the affinization $\Lgh$.
We end with some concluding remarks in  
Section \ref{ior}.  
In Examples \ref{E3},  \ref{E2} and  \ref{E4} we illustrate many of our results in the  case
where $(n,m)=(2,3)$.
\\ \\
While the main results of this paper are combinatorial, there is an easy 
consequence for representation theory.  We define a highest weight module $M$ to be {\it reflection complete} if  any sequence of adjacent Borel subalgebras induces proper inclusions 
between  submodules of $M$.  The submodules of $M$ obtained in this way are called {\it reflection submodules.} Let $\gr_0$  (resp. $\gr_1$) denote the half-sum of the positive even (resp. odd) roots of the distinguished Borel of $\fgo$ and set $\gr=\gr_0 -\gr_1$.  We show that 
the Vema module $M(-\gr)$ for $\fgo$ is reflection complete.  A related Verma module for the affinization at the critical  level is also reflection complete.  It follows that the poset of reflection submodules, the cover relation corresponds to an arrow in the graph $\cO^+(B,[{\stackrel{{\rm o}}{\fb}}{}^\dist,0])$.
\\ \\ We remark that groupoids are used to study Borel subalgebras of 
the affinization $\Lgh$ in 
\cite{GHS},10.2 and \cite{GK} using methods that are quite different from those used here.  I thank Maria Gorelik for some helpful correspondence.

\section{Young Diagrams, Groupoids and Borel subalgebras of $\fgl(n|m)$}\label{fc}
\subsection{The Lie superalgebra $\fgl(n|m)$}\label{Lsg}
Let $\fgo=\fgl(n|m)$.  The Cartan subalgebra $\fho$ of diagonal matrices in $\fgo$ is identified with $\ttk^{n|m}$.  
For $i\in [m+n], $  define $\gep_i\in \fhs$ so that $\gep_i(h) = h_i$ where  $h=\diag(h_1,\ldots,h_{m+n})\in \fho$.
For $i\in [n], $ we set $i' = i+m$ and $\gd_i = \gep_{i+m}.$ 
We have the usual permutation action  of $W$ on $\ttk^{n|m}$.  
The bilinear form $(\;,\;)$ is defined on $\fh^*$ is defined on the basis 
$\epsilon_{1},\ldots,   \epsilon_{n}, \delta_{1}, \ldots,\delta_{n}  $ by
\be \label{edform}(\epsilon_i,\epsilon_j) = \delta_{i,j} = - (\delta_i,\delta_j)\nn\ee 
and $(\epsilon_{i},\gd_{j})=0$ for all relevant indices $i,j$.
We write  $\Gl \in \ttk^{n|m}$ in the form
\be \label{ezr} \Gl = 
(a_1, \ldots ,a_n|b_1,\ldots,b_{m}) =  \sum_{i=1}^n  a_i
\gep_i - \sum_{j=1}^m  b_j \gd_j .\ee
If $\Gl$ is as above we have 
\be \label{exy} (\Gl,\gep_i- \gd_j)=  a_i-b_j. \ee 

\subsection{The Weyl groupoid}\label{aog}
In what follows, the details are slightly different from \cite{SV2} because we write all morphisms on the left of their arguments.  
Let  $W$ be a group of automorphisms of the groupoid  $\fG$ regarded as a category.
Then  $W$ acts on the base  $\fB$ of  $\fG$. We  define the  semi-direct product groupoid $W\ltimes \fG$  to have base $\fB$, and for $x,y \in \fB$, the set of morphisms from $x$  to $y$ is given by 
$$\mo(x,y) = \{(f, \gc)
\in W\ti \fG| f: \gc x \lra y \}.$$
We define the product
\be \label{Rq1}\mo(w,x) \ti \mo(x,y) \lra \mo(w,y)\ee
as follows.  Given $(f_1, \gc_1) \in \mo(x,y)$ and 
$(f_2, \gc_2)\in \mo(w, y)$, we have  morphisms  in $\fG$  
$\gc_1(f_2): \gc_1 \gc_2 w \lra \gc_1 x$ and $f_1:\gc_1 x \lra y.$
Thus the product  \eqref{Rq1} is defined by 
\be \label{Rq2}(f_1, \gc_1)\circ (f_2, \gc_2) =( f_1\circ \gc_1(f_2), \gc_1 \gc_2).\ee
We consider the Weyl groupoid $\cW$ constructed  
using the root system  $\Gd$. Let $\mathfrak T_{iso}$ be the groupoid with base 
$\Gd_{iso}$ the set of all the isotropic roots in $\Gd.$  
The non-identity  morphisms  are  $\gr_{\alpha}:\alpha \rightarrow -\ga$, $\ga\in\Gd_{iso}$ . The group $W$ acts on $\mathfrak T_{iso}$ in a natural way: $\alpha \rightarrow w(\alpha),\,
\gr_{\alpha} \rightarrow \gr_{w(\alpha)}$.  The {\it Weyl groupoid} $
\cW$ is defined
by 
\be \label{fgh}  \mathfrak{W} = W \coprod W \ltimes \mathfrak T_{iso},\nn\ee     the disjoint union of the group $W$ considered as a groupoid with a single point base $[W]$ and the semidirect product $W \ltimes \mathfrak T_{iso}$. 
Using \eqref{Rq2} it is easy to check that in $W \ltimes \mathfrak T_{iso}$ we have 
$(w,1)(1,\gr_\ga)(w^{-1},1) = (1, \gr_{w\ga})$, where 1 in the first and third factors stands for suitable identity mophisms.  We shorten this to 
\be \label{fah} w\gr_\ga w^{-1}=  \gr_{w\ga}.\ee 
Given a  functor $F:\mathfrak{W} \lra \mathfrak \fS(X)$,  we obtain
by restriction functors $F_1: W  \lra \mathfrak \fS(X)$ and 
$F_2: W \ltimes \mathfrak T_{iso} \lra \mathfrak \fS(X)$.  
If $X$ is a set with a natural $W$ action, $W\ti X\lra X, (w,x) = w\cdot x$, 
we assume that for all $w\in W, x \in X$, we have $F_1(w)x = w\cdot x$ and 
$F_2(w,1)x = w\cdot x$ whenever the latter is defined.  
Thus to define $F$ in these cases, we only need to specify $F(\ga)$ and $F(\gr_\ga)$ for 
$\ga\in \Gd_{iso}$ and check that the relation \eqref{fah} is preserved.

\subsection{Young diagrams and Borel subalgebras}\label{yd}  
 \subsubsection{Young diagrams, partitions and corners}\label{ypc}
Given positive integers $m, n$ with $m>n$, let $X$ be the set of Young diagrams that fit inside a $m\ti n$ rectangle {\bf R}.  Boxes in 
{\bf R} correspond to root spaces in  the 
upper triangular block $\fg^+_1$ of $\fgl(n|m)$.  Accordingly the rows of   
{\bf R} are indexed from the top down by $\gep_1,\ldots,\gep_n$  and columns of {\bf R} are indexed  from  left to right by  $\gd_1,\ldots,\gd_m$.  Thus the box in row $\gep_i$ and $\gd_j$
corresponds to the root space $\fg^{\gep_i-\gd_j}$ and we call this box $\fB(\gep_i-\gd_j)$, or often just $\gep_i-\gd_j$.  The set of roots of 
$\fg^+_1$ is $\Gd_1^+ = \{\gep_i-\gd_j|i\in [n], j \in [m]\}.$
\\ \\
A {\it partition} $\gl$ into $n$ {\it parts} is a  descending sequence 
\be \label{pdef}  
\gl=(\gl_1, \gl_2, \ldots \gl_n)\ee of non-negative integers.  
We call $\gl_i$ the $i^{th}$ {\it part} $\gl$. 
We also write 
\be \label{prf}\gl=(m^{a_m},  \ldots, i^{a_i}, \ldots)\ee
 to indicate that the part $i$ occurs with multiplicity ${a_i}$ in $\gl$.   Set $|\gl| = \sum_{i=1}^n \gl_i$. 
For a partition $\gl$,  the {\it dual partition} $\gl'$  has $j^{th}$ part 
$\gl'_j = |\{i|\gl_i \ge j\}|$. The Young diagram corresponding to $\gl$ has $\gl_i$ boxes in row $\gep_{n+1-i}$ with the first box on the left edge of {\bf R}. Our convention differs from the usual labelling of rows in a Young diagram, however we will also write matrices in {\bf R}, and it agrees with the usual way of counting rows in a matrix. To    
resolve this minor issue we  define the vector $\lst$ in\eqref{cd}, reversing the entries in $\gl$.
A partition is identified with the Young diagram it defines, but we say $\gl$ {\it is a partition} or $\gl$  {\it is a  Young diagram} according to which attributes of $\gl$ we wish to emphasize. 
\\ \\
If $\ga\in\Gd_1^+$,we say the box  $\fB=\fB( \ga)$ is an  {\it outer corner} of $\gl$ if $
  \fB \notin \gl$ 
and the addition of $ \fB$ to $\gl$ creates a new Young diagram.
Similarly box $ \fB$ is an  {\it inner corner} of $\gl$ if $
 \fB \in \gl$ 
and the removal of $ \fB$ from $\gl$ creates a new Young diagram.
Let 
$$X_\ga = \{\gl\in X| \fB(\ga) \mbox{ is an outer corner  of } \gl\}$$
and  
$$X_{-\ga} = \{\gl\in X| \fB(\ga) \mbox{ is an inner corner  of } \gl\}.$$  
For 
$\gl\in  X_{\ga}$,  $t_{ \ga}(\gl) \in X_{-\ga} $ is obtained by adding 
$ \fB(\ga)$ to $\gl$.  Similarly for 
$\gl\in  X_{-\ga}$,  $t_{- \ga}(\gl) \in X_{\ga} $ is obtained by deleting 
$ \fB(\ga)$ from $\gl$.  
This defines mutually inverse bijections  $ X_{\ga} \rl X_{-\ga} .$

\bl \label{rjs} 
If $\ga = \gep_i-\gd_{j}$, the following are equivalent
\bi \itema $ \fB(\ga)$ is an inner corner  of $\gl$, that is $\gl\in X_{-\ga}$.
\itemb $\gl_{n+1-i} = j$ 
and $\gl'_{j} = n+1-i$ .
\ei
If these conditions hold then 
\be \label{qef}  t_{-\ga}(\gl) = (\gl_1,\ldots , \gl_{n-i}, \gl_{n+1-i}-1,\gl_{n+2-i},\ldots, \gl_n).\ee
\el
\bpf  The equivalence of (a) and (b) follows since $ \fB(\ga)$ is an inner corner  of $\gl$ iff $ \fB(\ga)$ is the rightmost entry in row $\gep_i$ and the highest entry in column $\gd_{j}$.
The formula for $t_{-\ga}(\gl)$ follows from the definitions.\epf 
\noi Similarly we have
\bl \label{rjx} 
If $\ga = \gep_i-\gd_{j}$, the following are equivalent
\bi \itema $ \fB(\ga)$ is an outer corner  of $\gl$, that is $\gl\in X_{\ga}$.
\itemb $\gl_{n+1-i} = j-1$ and $\gl'_{j} = n-i$.
\ei
If these conditions hold then 
\be \label{qhf} t_{\ga}(\gl) = (\gl_1,\ldots , \gl_{n-i}, \gl_{n+1-i}+1,\gl_{n+2-i},\ldots, \gl_n).\ee
\el \noi 


 \subsubsection{Partitions, shuffles and paths}\label{bsa2}
Let 
\be\label{Imn}  I_0= \{1, \ldots, n\}, \quad I_1=\{1' , \ldots , m'  \}\nn\ee and 
let
$\cS_I$ be the group of permutations of $I= I_0 \cup I_1$.    
If $w \in \cS_I$ we define the {\it one line notation} for the permutation $w$ to be
\be \label{olnot}    {\it{{\underline{w}}}} = ( w(1), \ldots,w(n),w(1'),  \ldots w(m') ). \ee
The set of 
{\it shuffles} is
\[ \Sh = \{ \sigma \in \cS_I| 1, \ldots, n\;\; {\rm and} \;\; 1',\ldots,  m' \;\; {\rm are \; subsequences \;of}\;\;\underline{\sigma} \}. \]
We have a decomposition into right cosets
\be \label{drc}
\cS_I = W \ti \Sh.\ee 
\noi 
A {\it path} in {\bf R}  is a connected sequence of $m+n$ line segments of unit length starting at the 
top left corner and ending at the bottom right corner of
 {\bf R}.  Such  a path may be represented by a sequence  of symbols $\ur, \ud$ with 
$m$ $\ur$'s and $n$ $ \ud$'s. 
\bl \label{xmn} There is a bijection between $\gz:\Sh \lra X$ defined as follows: for a shuffle  $\gs$ we
draw a path in {\bf R}, where  the $k^{th}$ step is $\ud$ if $\gs(k) \le n$
and $\ur$ otherwise. Then the Young diagram $\gl= \gz(\gs)$ consists of the boxes below the path.
\el 
\bpf For (a) see  \cite{BN}  Lemma A.1   and Example A.2.\epf \noi 
The partition ${ \boldsymbol \gl}=
(m,1^{n-1})$ will play a special role. 
It is easy to see the following 
\bl \label{rags} \bi \itemo
\itema If $\gb= \gep_n-\gd_1$ is the unique simple root for the distinguished Borel, then 
${ \boldsymbol \gl}\subseteq \gl= \gz(\gs)$ iff $-\gb$ is the highest root of ${\stackrel{{\rm o}}{\fb}}{}(\gs)$. 
\itemb $\gl_1 =m$ iff 
$\gs(m')= n$.
\itemc
 ${ \boldsymbol \gl}\subseteq \gl$ iff $\gs(1)=1'$ and $\gs(m')= n$.
\ei
\el \noi


 \subsubsection{Borel subalgebras of $\fgl(n|m)$ and odd reflections}\label{0}

Let $\underline{{\stackrel{{\rm o}}{\cB}}}$ be the set of Borel subalgebras $\fbo$ of $\fgo=\fgl(n|m)$ 
containing the  subalgebra $\fho$ of diagonal matrices. 
 There is a bijection $\cS_I\lra  \underline{{\stackrel{{\rm o}}{\cB}}}$, 
$w\lra {\stackrel{{\rm o}}{\fb}}{}(w)$, where for $w\in \cS_I$  the corresponding Borel subalgebra 
${\stackrel{{\rm o}}{\fb}}{}(w) \in \underline{{\stackrel{{\rm o}}{\cB}}}$ has set of simple roots 
\be \label{pgs} 
{\stackrel{{\rm o}}{\Pi}}(w)=\{ \gep_{w(i)}- \gep_{w(i+1)}|i\in [m+n-1]\}.
\ee
Thus the distinguished Borel subalgebra ${\stackrel{{\rm o}}{\fb}}{}^\dist$
corresponds to $w=1.$  Now  we  have   bijections 
\be \label{bji} X \rl \Sh\rl  {\Fbo},\quad \quad \gl=\gz(\gs) \rl \gs\rl {\stackrel{{\rm o}}{\fb}}{}(\gs).\ee
\br \label{bjk} {\rm We expand on  another  description of the correspondence $\Fbo \lra X$ from  
\cite{M101} Section 3.4.   
Given ${\stackrel{{\rm o}}{\fb}}{}={\stackrel{{\rm o}}{\fb}}{}(\gs)\in{\stackrel{{\rm o}}{\cB}}$, write the
root $\ga = \pm(\epsilon_{i}-\delta_{j})$ in row $\gep_i$ and column $\gd_j$ of the rectangle {\bf R}  
where the sign is determined by the condition that $\fg^\ga \subset \fbo.$   
Next draw a path separating the roots of the form $\delta_{j} -\epsilon_{i}$ from those of the form $\epsilon_{i} - \delta_{j}.$  The
resulting array of roots together with the path is called the $\delta$-$\epsilon$-{\it diagram} for $\fbo$ (or for $\gs$). 
The path determines a Young diagram $\gl$ corresponding to $\fbo$. 
The simple roots of the form 
$\delta_{j} -\epsilon_{i}$ (resp.  $\epsilon_{i} - \delta_{j}$) are located in the inner (resp. outer) corners of $\gl.$  We also call these corners 
inner and outer corners of the path. By the shuffle condition, between an inner corner and the next outer corner, we find only simple roots of the form $\epsilon_{\ell} - \epsilon_{\ell+1}$ and between an outer corner and the next inner corner, only those of the form $\delta_{\ell}-\delta_{\ell+1}.$
} 
\er \noi 
If $\ga = \gep_i-\gd_{j}$, we determine the analogs for $\Sh$ and ${{\stackrel{{\rm o}}{\cB}}}$ of the bijections 
$t_{\pm\ga}: X_{\pm\ga} \lra X_{\mp\ga} $ from Subsection \ref{ypc}. 
Define \be \label{323} 
\Sh(\ga)= \{\gs\in \Sh|\ga \in{\stackrel{{\rm o}}{\Pi}}(\gs)\},\quad {\stackrel{{\rm o}}{\cB}}(\ga)=  \{ \fbo(\gs)| \gs \in \Sh, \;\ga\in 
{\stackrel{{\rm o}}{\Pi}}(\gs)\}.\ee
Let $\fbo\in {\stackrel{{\rm o}}{\cB}}(\ga)$ be a Borel with simple roots ${\stackrel{{\rm o}}{\Pi}}$. 
For $\gb \in \pigo$ we define a root $r_\ga(\gb)$ by
\[ r_\ga(\gb) = \left\{ \begin{array}
  {rcl}
  -\ga  & \mbox{if} & \gb = \ga \\
\ga + \gb & \mbox{if} & \ga + \gb \mbox{ is a root}.\\
\gb & \quad \quad \mbox{otherwise} &
\end{array} \right. \]
Then $r_\ga(\pigo)$ is the set of simple roots for a Borel subalgebra $r_\ga(\fbo)$. Clearly
$r_{\pm\ga}:{\stackrel{{\rm o}}{\cB}}(\ga)\lra {\stackrel{{\rm o}}{\cB}}(-\ga)$ are mutually inverse bijections.  
If  $\ga =\gep_i-\gd_j$, it follows from \eqref{pgs} and \eqref{323} that 
$\gs\in \Sh(\ga)$ (resp. $\gs\in \Sh(-\ga)$)  iff $i$ is the predecessor (resp. successor) of  $j'$in $\underline{\sigma}$. Define inverse bijections $r_{\pm\ga}:\Sh(\pm\ga)\lra \Sh(\mp\ga)$ by the requirement that if $\gt = r_{\pm\ga}(\gs)$ then the order of $i, j'$ in $\underline{\sigma}$ 
and $\underline{\gt}$ is interchanged. As a product of permutations we have
\be \label{32r} r_{\pm\ga}(\gs) = (i,j')\gs.\ee  We say that 
$r_\ga(\fbo)$ 
 (resp.  $r_{\ga}(\gs)$) is obtained from $ \fbo$ (resp. $\gs$) using the {\it odd reflection} $r_{\ga}$. As in \cite{M101} Chapter 3, we also say that  $ \fbo$  and $r_\ga(\fbo)$ are {\it adjacent} Borel subalgebras.

\bl \label{32s} 
The maps from \eqref{bji} restrict to bijections 
\be \label{bj1} X_{\pm\ga} \rl \Sh(\pm\ga)\rl  {{\stackrel{{\rm o}}{\cB}}}(\pm\ga),\nn\ee
and we have a commutative diagram.
\[
\xymatrix@C=2pc@R=1pc{
X_{\ga}\ar@{<->}[rr] \ar@{<->}_{t_\ga} [dd]&&
\Sh(\ga)\ar@{<->}_{r_\ga}[dd]&&
{\stackrel{{\rm o}}{\cB}}(\ga) \ar@{<->} [ll] \ar@{<->}^{r_\ga}[dd]&\\ \\
X_{-\ga}\ar@{<->}[rr] &&\Sh({-\ga})\ar@{}[rr] &&
{\stackrel{{\rm o}}{\cB}}({-\ga})\ar@{<->}[ll] &}
\] 
\el
\bpf The correspondence $\Sh(\ga)\rl  {{\stackrel{{\rm o}}{\cB}}}(\ga)$ and the commutativity of the right square are  immediate from
\eqref{323}. 
Consider the diagram below. 
On the left (resp.  right) we show the paths determined  $\gl= \gz(\gs)$  
when 
$\ga = \gep_i-\gd_{j}$ is an outer (resp. inner) corner of  $\gl= \gz(\gs)$.  

\Bc
\setlength{\unitlength}{0.8cm}
\begin{picture}(10,3)(-1,-1.50)
\thinlines
\linethickness{0.05mm}
\put(1.39,0.35){$\ga$}
\multiput(1.0,0)(1,0){2}{\line(0,1){1.0}}
 \multiput(1,0)(0,1){2}{\line(1,0){1.0}}
  \linethickness{0.5mm}
\put(8.0,1){\line(1,0){1}}
\put(9.0,1){\line(0,-1){1}}
\put(6.45,0.35){$ \gep_i$}
\put(7.19,0.35){\dots}
\put(8.38,-.75){\vdots}
\put(0.19,0.35){\dots}
\put(1.38,-.75){\vdots}
\put(-.45,0.35){$ \gep_i$}
\put(8.39,-1.35){$ \gd_j$}
\put(1.39,-1.35){$ \gd_j$}
  \linethickness{0.05mm}
\put(8.39,0.35){$\ga$}
\multiput(8.0,0)(1,0){2}{\line(0,1){1.0}}
 \multiput(8,0)(0,1){2}{\line(1,0){1.0}}
  \linethickness{0.5mm}
\put(1.0,0){\line(1,0){1}}
\put(1.0,1){\line(0,-1){1}}
  \end{picture}
\Ec
The remaining statements  follow from Lemma \ref{xmn} and a consideration of these diagrams.
\epf

 \subsubsection{Isomorphisms of groupoids}\label{iog}
 Now let $\fT(X)$ be the subgroupoid of $\mathfrak \fS(X)$ with base $\{X_{\pm\ga}|\ga\in\Gd_1^+\}$ and non-identity morphisms $\{t_{\pm\ga}|\ga\in\Gd_1^+\}$. 
There is a 
functor $F:\fT_{\iso}\lra \fT(X)$ as in \eqref{cwr} given by 
$ F({\alpha})=X_\ga$  and 
$F(\gr_{\alpha}) =t_\ga$ for  ${\alpha}\in \pm \Gd_1^+$. 
If $\gl \in X_{\pm\ga}$ we have
\be \label{xit}  F(\gr_{\pm\ga})(\gl)=t_{\pm \ga}(\gl).\ee
Since all we have done is rename the objects and morphisms, it is clear that  $F$ is an isomorphism of groupoids.\\ \\  Similarly
let $\fT({\stackrel{{\rm o}}{\cB}})$  be the groupoid with base $\{{\stackrel{{\rm o}}{\cB}}(\pm\ga)|\ga\in\Gd_1^+\}$ and non-identity morphisms $\{r_{\pm\ga}|\ga\in\Gd_1^+\}$. 
There is a 
functor $B:\fT_{\iso}\lra \fT({\stackrel{{\rm o}}{\cB}})$ given by 
$ B({\alpha})={\stackrel{{\rm o}}{\cB}}(\ga)$  and 
$B(\gr_{\alpha}) =r_\ga$ for  ${\alpha}\in \pm\Gd_1^+$. We have a commutative diagram of isomorphisms.
\be \label{bby}
\xymatrix@C=2pc@R=1pc{
\mathfrak  T_{iso}\ar@{<->}_= [dd]&&
\fT(X)\ar@{<-}_F[ll]: \ar@{<->}^{\cong}[dd]&\\ \\
\mathfrak  T_{iso} \ar@{-}[rr]_B &&
\fT({\stackrel{{\rm o}}{\cB}})\ar@{<-}[ll] &}
\ee
\section{Extending the action of $\mathfrak  T_{iso}$} \label{APG}
\subsection{A general construction} \label{AP1}
Let $\cH$ be  the smallest groupoid that is not a disjoint union of groups.  If there is  an action of the disjoint union $\cH \coprod \mathfrak T_{iso}$ on a set $\X$ 
satisfying Hypothesis \ref{ky} below, we show the action of 
 $\mathfrak T_{iso}$ can be extended to an action on
the equivalence classes  of $\X\ti\Z$ under a certain equivalence relation.  
\\ \\
Let $\nu$ be the permutation of $\{1, \ldots, n\}$ given by 
$\nu(k) = k-1 \mbox{ mod } n$. Recall notation \eqref{ezr}. 
We extend $\nu$ to  an operator on $\ttk^{n|m}$  by setting
\be \label{e67l}\nu(a_1, \ldots, a_n|b_1, \ldots, b_m) = (a_{\nu(1)}, \ldots, a_{\nu(n)}|b_1, \ldots, b_m) \ee
Note that $\nu$ preserves the bilinear form defined in Subsection \ref{Lsg}: we have $(\nu\Gl,\nu\ga)=(\Gl,\ga)$ for  $\Gl, \ga \in \ttk^{n|m} $. 
Also \eqref{e67l} implies 
\be \label{e67m}\nu(\gep_i - \gd_j) = \gep_{i+1}-\gd_j  \mbox{ for } i\in [n-1], \quad \nu(\gep_n- \gd_j)= \gep_1 - \gd_j.\nn\ee
Denote the objects of $\cH$ 
by $\da, \pa$ and the non-identity mophisms by
  ${}^{-}:\da\lra\pa $ and   ${}_{-}:  \pa\lra\da$.  
Suppose there is a 
 functor $\F: \cH \coprod \mathfrak T_{iso}\lra \fS(\X)$.  Define $\X_\gb=\F(\gr_\gb)$ for 
$\gb\in \Gd_\iso$, $\X(\da)=\F(\da)$ and $\X(\pa)=\F(\pa) $.
Set 
  $\overline{\mu} =\F(^{-})(\mu) $ and   $\underline{\mu} =\F({}_{-}) \mu$
for 
$\mu\in  \X(\da)$, $\mu\in \X(\pa)$ respectively.
\bh \label{ky} 
  Assume that   the following condition holds.
\bi \itemo
If $\gb\in \Gd_{1}^+$, $\gl\in \X_{\gb}\cap  \X(\da)$ and 
$\overline{\gl} \in \X_{\nu\gb}$, then ${\F(\gr_{\gb})(\gl)}\in  \X(\da)$ and   
\be \label{kz} \overline{\F(\gr_{\gb})(\gl)} = \F(\gr_{\nu\gb})(\overline{\gl}).\ee \ei\eh
\bl If Hypothesis \ref{ky} holds, then 
if  $\gb\in \Gd_{1}^+$,  $\mu\in  \X_{-\nu\gb}\;\cap \in \X(\pa)$ and  $\underline{\mu} \in \X_{-\gb}$, then  ${\F(\gr_{-\nu\gb})(\mu)}\in  \X(\pa)$ and  
\be \label{kf} \underline{\F(\gr_{-\nu\gb})(\mu)} = \F(\gr_{-\gb})(\underline{\mu}).\ee
\el 
\bpf Set $\mu = \overline{\F(\gr_{\gb})(\gl)}= \F(\gr_{\nu\gb})(\overline{\gl})$ and solve  
 for $\gl$ in two different ways.
\epf
\noi 
Let $\sim$ be the smallest equivalence relation on 
$\X\ti \Z$ such that $(\gl,k)\sim (\overline \gl,k+1)$ if $\gl_1=m$.   
Denote the equivalence class of $(\gl,k)\in \X\ti Z$ under $\sim$ by $[\gl,k]$ and let  
$[\X\ti Z]$ be the set of equivalence classes.  In the pair $({\gl},k)  \in X\ti \Z$, 
 we call $k$ the {\it rotation number}.
To extend the action of $\mathfrak  T_{iso}$ on $X$ given by \eqref{xit} 
to an action on $[\X\ti Z]$, we define a functor 
\be \label{krt} F:\mathfrak  T_{iso}\lra \mathfrak S[\X\ti \Z].\ee
In the above $\fS[\X\ti \Z]$ means $\fS([\X\ti \Z])$.  We make similar abbreviations without further comment. 

\noi 
For $\ga\in \Gd_\iso$   set
\be \label{sox}Y_\ga= \bcu_{j\in \Z} [\X_{ \nu^{j}\ga},j ].\ee
\bt \label{shg}  If the functor $\F$ satisfies Hypothesis \ref{ky}, then 
\bi \itema There is a functor as in \eqref{krt} given by 
\be \label{syx}F(\ga) = Y_\ga
\mbox{ and } F(\gr_\ga)[\gl, j] = [\F_{ \nu^{j}\ga}(\gl),j)]
\mbox{ for  } (\gl, j)\in (X_{ \nu^{j}\ga},j )
.\ee
\itemb 
The injective map 
$\X\lra [\X\ti Z],$ $\gl\lra [\gl,0]$ is compatible with the 
 $\mathfrak T_{iso}$-action.  If 
$\gl \in X_{\ga}$, then $F(\gr_{\ga})[\gl,0]= [\F(\gr_{\ga})(\gl),0].$
\ei
Thus the action of 
 $\mathfrak T_{iso}$ on $\X$ extends to an action on $[\X\ti Z].$
\et
\bpf (a) Since $(\gl,j)\sim  (\overline{\gl},j+1)$, to  show the action in \eqref{syx} is well-defined, we need to show the following.
If 
$\gl \in X_{ \nu^{j}\ga}$ and  $\overline{\gl}\in X_{ \nu^{j+1}\ga}$ then
$$(\F(\gr_{ \nu^{j}\ga})(\gl),j)\sim 
(\F(\gr_{ \nu^{j+1}\ga})(\overline{\gl}),j+1).$$
If  $\gb={ \nu^{j}\ga}$, then
using the definition of $\sim$ 
 and \eqref{kz}
\by \label{f2}   
(\F(\gr_{\gb})(\gl),j)
&\sim& 
(\overline{\F(\gr_{\gb})(\gl)},j+1)
\nn\\&=& 
(\F(\gr_{\nu\gb})(\overline{\gl}),j+1).
\nn
\ey  as required.  
Statement (b)  is \eqref{syx} with $j=0$.\epf

\subsection{Application to Young diagrams and Borel subalgebras} \label{AP2} As before $X$ is the set of Young diagrams that fit inside the  $m\ti n$ rectangle {\bf R}.
Set \be\label{yyt}X(\da) = \{\gl\in X|\gl_1=m\}, \quad X(\pa) = \{\gl\in X|\gl'_1<n\}.\ee 
If $\gl \in X(\da)$, 
let   $\overline \gl =  (\gl_2, \ldots,\gl_n,0 )$  be the partition obtained from removing the first part from $\gl$  
and if $\gl \in X(\pa)$, set 
$\underline \gl =(m,\gl_1,\gl_2,\ldots ,\gl_{n-1}).$  We have $\overline{\gl}_ {\nu(k)}={\gl}_k$ if $k\in [n-1]$. 
As a Young diagram 
$\overline \gl $
is obtained from  $\gl$ by deleting the bottom row.   
Passing from $\gl$ to 
$\overline \gl $  allows us to continue the process of forming larger Young diagrams by adding outer corners.  
\\ \\
Define the  functor 
$\F:\cH \coprod \mathfrak T_{iso}\lra \fS(X)$ 
by the requirements that  
\be\label{yit}
\F(\da) =  X(\da) , \quad 
\F(\pa) =X(\pa), \quad \F({}^{-})(\gl) =\overline \gl, \quad
\F({}_{-}) (\gl)=\underline \gl\ee
and the restriction $F$ of $\F$ to 
$\mathfrak T_{iso}$ is given by \eqref{xit}.
\bl \label{rhy} The functor $\F$ satisfies Hypothesis \ref{ky}.
\el
\bpf Suppose $\gb= \gep_i-\gd_{j}\in \Gd_{1}^+$.  
We have to show that if 
$\gl\in X_{\gb}\cap X(\da)$ and $\overline{\gl} \in 
X_{\nu\gb}$, then ${t_{\gb}(\gl)}\in X(\da)$ and  
\be \label{kzt} \overline{t_{\gb}(\gl)} = t_{\nu\gb}(\overline{\gl}).\ee
Since $\gl_1 = m$, $\gl$ has no outer corners in row $\gep_{n}$. Thus 
 $i<n$.
The equivalence of (a) and (b) in  Lemma \ref{rjx} implies 
$\gep_{i+1}-\gd_j$  is an outer corner of $\overline{\gl}$ and the rest follows from \eqref{qhf}.  
\epf 
\noi Recall the bijection 
$\gz:\Sh \lra X$ given in Lemma \ref{xmn}. 
 Suppose  $\gl=\gz(\gs)$ is a partition with $\gl_1 = m$.  By Lemma \ref{rags} 
  this means that  $\gs(m+n) = n$.  
We determine the shuffle that corresponds to  
$\overline{\gl}$. 
Extend the    permutation    $\nu $   of $I_0$ from Subsection \ref{AP1}, to a  permutation of $ I $ by setting $\nu j =j$, $j\in I_1.$
 Let  $\overline{\sigma}$ be the   permutation whose one-line notation is
$$\overline{\underline{\sigma}}=(1,  \nu^{-1}\gs(1), \ldots,\nu^{-1}\gs(n),\nu^{-1}\gs(1'),  \ldots \nu^{-1}\gs(m+n-1) ).$$

\bl \label{ymn} We have 
\bi \itema
$\overline{\gl} =\gz(\overline{\sigma}).$
\itemb  Suppose $\ga =\gep_i-\gd_j$ and  
$r_{\ga}(\gs) = (i,j')\gs$ is obtained from $\gs$ using the { odd reflection} $r_{\ga}$.
Then  
\be \label{e6z}\overline{r_\ga\gs}= r_{\nu\ga}\overline{\gs}.\ee 
\ei\el
\bpf   Since the first entry in $\overline{\sigma}$ is 1, the first step in the path corresponding
to $\gz(\overline{\sigma})$ is $\ud$.  
For $i\ge 0$, step $i+1$ of this path is the same as 
step $i$ of the path for $\gl$.  This proves (a).  For (b), suppose $\gs(k)= i, \gs(k+1)= j'$, set 
$\gs'=r_{\ga}(\gs)$ and 
consider the permutations given in one-line notation by 
$$(n, \gs(1), \ldots, \gs(k-1), i, j',\gs(k+2), \ldots)$$
$$(n, \gs(1), \ldots, \gs(k-1), j',i,\gs(k+2), \ldots)$$
To obtain $\overline{\sigma}$  and $\overline{\sigma'}$ from these, all unprimed entries are decreased by 1 mod $n$.  Equation \eqref{e6z} follows from this.
\epf
\noi 
As in Subsection \ref{0}, let ${\stackrel{{\rm o}}{\cB}}$ the set of Borel subalgebras in $\fgo$
having the upper triangular matrices as the even part.  
Set
$$ \cB(\da) = \{\fbo(\gs)| \gl=\gz(\gs)\in X(\da)\}, \quad 
\cB(\pa) = \{\fbo(\gs)| \gl=\gz(\gs)\in X(\pa)\}$$
We extend the functor $B:\fT_{\iso}\lra \fT(\Fbo)$ from Subsection \ref{iog} to a functor 
$\B:\cH \coprod \mathfrak T_{iso}\lra \fS(\Fbo)$ 
such  that  
\be\label{yat}
\B(\da) =\cB(\da), \quad 
\B(\pa) = \B(\pa), \quad \B({}^{-})(\fbo(\gs)) =\fbo(\overline \gs), \quad
\B({}_{-}) (\fbo(\gs)) =\fbo(\underline \gs)(\gl),\ee
where we  use the bijections $\gl=\gz(\gs) \rl \gs\rl \fbo(\gs)$ from \eqref{bji} and 
 $\overline \gl=\gz(\overline \gs),$  $\underline \gl=\gz(\underline\gs) .$ 
\bl \label{rhk} The functor $\B$ satisfies Hypothesis \ref{ky}.
We have 

\be \label{kx} \overline{\B(\gr_{\ga})(\fbo(\gs))} = \B(\gr_{\nu\ga})(\fbo(\overline{\gs})).\ee\el
\bpf By Lemma \ref{32s} we can rewrite  \eqref{e6z} in the form 
\be \label{kxz}
\overline{r_\ga(\fbo(\gs))}= r_{\nu\ga}(\fbo(\overline{\gs})).\ee
This is equivalent to \eqref{kx}.
\epf

\bc \label{sz} \bi \itema
There is a functor $F:\fT_{\iso}\lra \fS[X\ti \Z]$ such that for $\ga\in \Gd_{\iso}$
$$F(\ga) = \bcu_{j\in \Z} [X_{ \nu^{j}\ga},j ]
\mbox{ and } 
F(\gr_\ga)[\gl, j] = [t_{ \nu^{j}\ga}(\gl),j)]$$
for $(\gl, j) \in (X_{ \nu^{j}\ga},j).$
\itemb
There is a functor $B:\fT_{\iso}\lra \fS[\Fbo\ti \Z]$ such that for $\ga\in \Gd_{\iso}$
$$B(\ga) = \bcu_{j\in \Z} [\Fbo({ \nu^{j}\ga}),j ]
\mbox{ and } 
B(\gr_\ga)[\gl, j] = [r_{ \nu^{j}\ga}(\gl),j)]$$
for $(\gl, j) \in ({\stackrel{{\rm o}}{\cB}}({ \nu^{j}\ga}),j).$
\itemc The functors $F, B$ that fit into the commutative diagram \eqref{sez}
\ei
\ec
\bpf  This follows directly from Lemmas \ref{rhy}, \ref{rhk} and Theorem  \ref{shg}.
\epf 

 \subsection{Pseudo-corners}\label{gbx}
Let 
$$X^{\red}= \{\gl \subseteq {\bf R} | \gl_1 <m\} = X \bsk   X(\da).$$ This is the set of {\it reduced diagrams}, the set of 
 {\it strongly reduced diagrams} is, 
$$X^{\str}= \{\gl \subseteq {\bf R} |  \gl_n = \gl'_m =0\}.$$ 
If ${ \boldsymbol \gl} \subseteq \gl$,
we call $\gep_n-\gd_1$   a   {\it outer pseudo-corner} of  $\gl$, and if 
$\gl \in X^{\str}$ we call $ \gep_1-\gd_m  $   a {\it  inner  pseudo-corner} of  $\gl$.  These are not outer  or inner corners, but behave in some ways as if they were.  Pseudo-corners account for the cover relations in the poset $[X\ti \Z]$ that do not arise from $X$  and Theorem \ref{shg} (b).  They also  play an important role in the next Section, see 
Lemmas \ref{c1} and \ref{c2}. It is easy to see the following.
\bl \label{rey} \bi \itemo
\itema  If 
$\ga =\gep_n-\gd_1$ 
and ${ \boldsymbol \gl} \subseteq  \gl$,
 then $\ga$ is not an outer corner of $\gl$, but 
$\nu\ga = \gep_1-\gd_1$ is an outer corner of 
$\overline{\gl}.$
\itemb  
If $\ga =\gep_1-\gd_m$ and $\gl \in X^{\str}$, then $\ga$ is not an inner corner of $\gl$, but $\nu^{-1}\ga=\gep_n 
-\gd_m$ is  an inner corner of $\underline{\gl}.$
\ei
\el  \noi
We pay special attention to pseudo-corners in Corollary \ref{sz}.  Suppose 
$\ga =\gep_n-\gd_1$ then  
$\gb:=\nu\ga = \gep_1-\gd_1$ is an outer corner of 
$\overline{\gl}.$   We have $(\gl,0)\sim(\overline{\gl},1).$ Hence 
\be \label{1no} F(\gr_\ga)[\gl,0]=F(\gr_\ga)[\overline{\gl},1]= [t_\gb\overline{\gl},1].\ee
Now suppose that  $\ga =\gep_1-\gd_m$ and $\gl \in X^{\str}$, then $\ga$ is not an inner corner of $\gl$, but $\gb=\nu^{-1}\ga=\gep_n 
-\gd_m$ is  an inner corner of $\underline{\gl}.$  We have 
$(\gl,0)\sim(\underline{\gl},-1)$ and 
\be \label{2no}F(\gr_{-\ga})[\gl,0]=F(\gr_{-\ga})[\underline{\gl},-1]= [t_{-\gb}\underline{\gl},-1].\ee

 \subsection{$[X\ti \Z]$ as a graded poset}\label{gax}
We give the pair $({\gl},k)$ a degree by setting
\be \label{e67n}\deg({\gl},k)=|\gl| + km\in \Z.\ee  This degree passes to equivalence classes.

\bl \label{rjg1} 
$\deg(\overline{\gl},k+1)=\deg(\gl,k)$ if $\gl_1=m$. 
\el  \bpf This is clear since $\gl$ contains $m$ more boxes than $\overline{\gl}$. 
\epf  \noi Thus 
$[X\ti \Z]= \bcu_{r\in \Z} [X\ti \Z]_{r}$ is a graded poset where $[X\ti \Z]_r = \{[{\gl},k]|\deg[{\gl},k] =r\}.$ It is clear that $[X\ti \Z]$ enjoys a periodicity property: there is a bijection $[X\ti \Z]_{r}\rl [X\ti \Z]_{r+mn}$,  $[{\gl},k]\rl [{\gl},k+n]$. 
 \bexa \label{E3}{\rm  
 Let $(n,m) =(2,3)$.  Below we give the portion of the graph/graded poset in $[X\ti Z]$ degrees 0-6.  The degree, as defined in \eqref{e67n}  is given in the last row. Because of the periodic nature of this poset, the part in degree $i+6$ is obtained  from the degree $i$ part by increasing the rotation number by 2. So the portion below determines the entire graph.   
 In examples, the rotation number can be hidden by adding the superscript $+$ (resp. $-$)  whenever the roation number is increased (resp. dereased) by one.   Assume below that unadorned Young diagrams have rotation number 0.  When for example  $\gl=(2,1)$, we use  $\gl^-$ and $\gl^+$ as shorthand for $((2,1),-1)$ and $((2,1),1)$ respectively.  
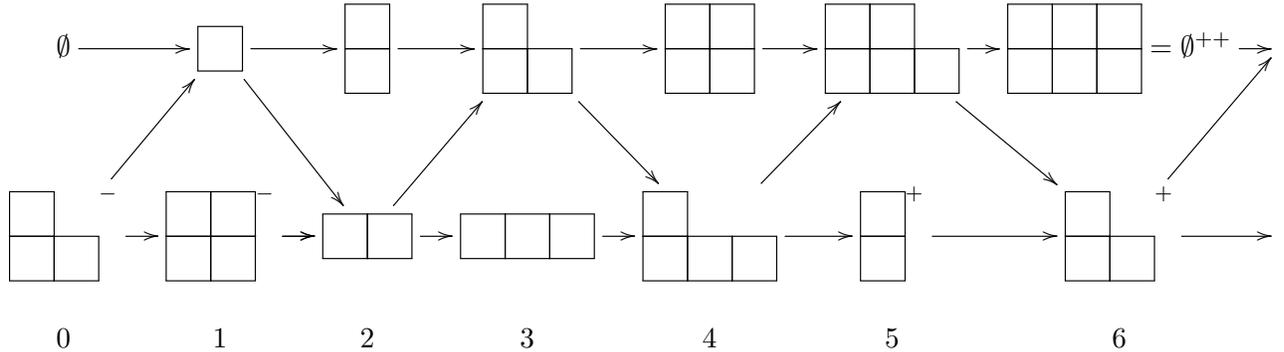
\begin{equation*}
\xymatrix@C=1pc@R=1pc{\emptyset\ar@{->}[r] &
\ydiagram{1}\ar@{->}[ddr]&
\ydiagram{1,1} \ar@{<-}[l] \ar@{->}[r]& \ydiagram{0+1,2}\ar@{->}[r] \ar@{->}[ddr]   &
\ydiagram{2,2} \ar@{->}[r]&\ydiagram{0+2,3} \ar@{->}[r] \ar@{->}[ddr]& \ydiagram{3,3}=\emptyset^{++}\ar@{->}[r]&
\\ \\
 \ydiagram{0+1,2}^{-}\ar@{->}[r] \ar@{->}[uur]&\ydiagram{2,2}^{-}\ar@{->}[r]&
\ydiagram{2}\ar@{<-}[l] \ar@{->}[r] \ar@{->}[uur]&\ydiagram{3}\ar@{->}[r]&
\ydiagram{0+1,3}\ar@{->}[r] \ar@{->}[uur]&\ydiagram{1,1}^+\ar@{->}[r]& \ydiagram{0+1,2}^+\ar@{->}[r]\ar@{->}[uur]&\\
0&1&2&3&4&5&6}
\end{equation*}
 The most interesting behavior occurs when the rotation number changes after a morphism is applied.  Here are some examples 
By  \eqref{1no} 
we obtain 
$$F(\gr_{\gep_2 -\gd_1})[(3,i),0)]=[t_{\gep_1 -\gd_1}(i,0),1] =[(i,1)^{+}] \mbox{ for } i=1,2.$$   
Together with the more obvious morphism $(1,1)^+ \lra (2,1)^+$, this accounts for the extra arrows in the lower right that do not appear in the Hasse diagram for the poset $X$.
Using   \eqref{2no} 
$$F(\gr_{\gd_3-\gep_1})[(i,0),0)]=F(\gr_{\gd_3-\gep_1})[(3,i),-1)]=[t_{\gd_3-\gep_2 }(3,i),-1]  =[(2,i)^-]  \mbox{ for } i=1,2.$$
The diagram shows only positive morphisms, but the last equation is equivalent to
$$F(\gr_{\gep_1-\gd_3})[(2,i)^-] =[(i,0),0)] \mbox{ for } i=1,2.$$
This accounts for the  extra arrows in the lower left  of the diagram.
}\eexa


\section{The Sergeev-Veselov Functor }\label{affc}

\subsection{The action on $\ttk^{n|m}$}\label{aog11}
In  \cite{SV101} Sergeev and Veselov defined 
an action of the Weyl groupoid  $\mathfrak{W}$  on $\ttk^{n|m}$ 
depending on a non-zero parameter $\gk$. 
If  $\gk = -1$ all orbits are described in \cite{M22}  Theorem 5.7. 
By  \cite{SV101} Theorem 3.2,  $\cW$ has an infinite orbit  if (and only if) $\gk$ is rational with $-1\le \gk <0.$  
The proof of this statement reduces to the case $\gk=-n/m$.  We consider only the case
$\gk=-n/m$ and  
our notation  is slightly different from \cite{SV101}. It will allow us to consider vectors and matrices with integer entries. Starting in Subsection \ref{cpc} we assume that $m,n$ are coprime and $m>n.$ We define the {\it Sergeev-Veselov Functor} 
\be \label{cwr1}SV:\mathfrak{W}\lra \mathfrak \fS(\ttk^{n|m}).\ee
With  the usual permutation action 
 of $W$ on $\ttk^{n|m}$,  it suffices to describe 
the restriction of $SV$ to $\mathfrak  T_{iso}$.  
For $\ga=\gep_i- \gd_j \in \Gd_\iso^+$, 
 define $v_\ga = n\gep_i -m \gd_j$ and 
\be \label{e99}\Pi_{\ga} = \{\Gl\in \ttk^{n|m}| (\Gl, \ga)=0\},\quad \Pi_{-\ga} = \{\Gl\in \ttk^{n|m}| (\Gl,\ga)=n-m\}.\ee
Next define $\gt_{\pm\ga} :\Pi_{\pm\ga} \lra \Pi_{\mp\ga} $ by 
\be \label{e9}\gt_{\pm\ga} (\Gl) = \Gl \pm v_\ga.\ee 
The 
functor  $SV$ from \eqref{cwr1} is given  by 
$SV(\pm \ga) =\Pi_{\pm\ga}$, and $SV(\gr_{\pm\ga}) =\gt_{\pm\ga}$. 
This action restricts to an action of $\cW$ on $\Z^{n|m}$. Abusing notation,  the restriction of $SV$ to $\mathfrak  T_{iso}$ gives a functor 
\be \label{e11}SV:\mathfrak  T_{iso}\lra \mathfrak \fS(\Z^{n|m}).\ee 
If $\ga\in \Gd_\iso^+$, we call $\gt_{\ga}$ (resp. $\gt_{-\ga}$) a positive (resp. negative) morphism.
\\ \\
For $\Gl$ as in \eqref{ezr}, set
\be \label{e73}
^L\Gl = {(a_1,\ldots a_{n})}^\trans  , \quad ^R\Gl = ( b_1,\ldots ,b_m)
.\nn\ee 
and 
\be \label{e83}
|^L\Gl| = \sum_{i=1}^n a_i, \quad |^R\Gl| = \sum_{i=1}^m b_i.\nn\ee 

\bl \label{fj} If $\Gl$  has integer entries and the entries  in  
$^L\Gl$ $($resp. $^R\Gl$$)$ form a complete set of representatives mod $n$ $($resp. mod $m$$)$, then the same is true for the entire orbit of $\Gl$.
\el
\bpf The condition is invariant under permutations of the entries in $\Gl$, so the result holds by \eqref{e9}.
\epf \noi 
We consider the $\mathfrak  T_{iso}$-orbit 
$\cO$ 
of 
   \be \label{e7}\Gl_0 = (m(n-1), \ldots,m,0| 0,n,\ldots,n(m-1))
= \sum_{i=1}^n  m(n-i)
\gep_i - \sum_{j=1}^m n(j-1)  \gd_j. 
\ee under the action from \eqref{e11}. 
 \subsection{The augmented matrix $\widehat\A(\Gl)$}\label{fi}  
Let $\A(\Gl)$ be the matrix with entry in row $\gep_i$ and column $\gd_j$ given by 
\be \label{Adef}
\A(\Gl)_{i,j}= (\Gl,\gep_i-\gd_j) =a_i -b_j.\nn
\ee 
For $a\in\ttk$, set 
$\bA=(a,\ldots, a|a,\ldots, a) = a \;\str$ where $\str$ is the supertrace. 
\bl \label{l73}
Given $\Gl, \Gl'\in \ttk^{n|m}$, we have $\A(\Gl)=\A(\Gl') $ iff   $\Gl=\Gl' +\bA$ for some $a\in \ttk$.\el 
\bpf  If $\A(\Gl)=\A(\Gl') $, then $a= a_i -a'_i  =  b_j -b'_j$ is independent of $i,j$. 
\epf
\noi Since some information is lost in passing from $\Gl$ to $\A(\Gl)$,  we also use the (doubly) augmented matrix. 

\[\widehat \A(\Gl) =\begin{tabular}{|c||ccc|} \hline
&  &$^R\Gl$&
\\ \hline\hline &     &                          & \\ $^L\Gl$
&&$\A(\Gl)$
 & \\
 & &                              & \\
  \hline
\end{tabular}\]
In Example \ref{E2}, when $(n,m) =(2,3)$, we give several matrices $\widehat \A(\Gl)$. %
\subsection{Construction of $\Gl=x(\gl) \in \Z^{n|m}$ from $\Gl_0$ and $\gl \subseteq$ {\bf R}} \label{coL}
We use 
$\Gl_0$ and $\gl \subseteq$ {\bf R}, to construct certain elements 
$x(\gl)\in \cO$.  However we point out that instead of $\Gl_0$, we could have used  
   \be \label{e7a}\Gl_0^{(k)} = \Gl_0 + k \bf{m}
\ee 
as a starting point for any $k\in \Z$.  Passing from $\Gl_0$ to $\Gl_0^{(k)} $,  the entries of $\Gl_0$ all increase by the same integer, so the difference between entries is unchanged. Thus $\A(\Gl_0^{(k)}) = \A(\Gl_0)$ and the key Equations \eqref{e67d} and \eqref{e67c} still hold when 
$ \Gl_0 $ is replaced by $\Gl_0^{(k)}$. In addition $\Gl_0^{(kn)}\in \cO$ and  $\Gl_0^{(k)}\in \cO$ up to rotation  by Corollary \ref{aed} and  \eqref{xdef}.
 \\ \\
Given $\gl \subseteq$ {\bf R}, consider
\be \label{cd}  
\gl'=(\gl'_1,\ldots, \gl'_m) 
\mbox{ and } \lst =(\gl_n,\ldots, \gl_1)^\trans\ee 
as vectors.  Reversing the entries in $\gl$ means that when written as a column vector alongside $\gl \subseteq$ {\bf R}, the $i^{th}$ entry of $\lst$ gives the number of entries, $\gl_{n+1-i}$
in row $\gep_i $ of $\gl$. Define $\Gl =x(\gl)\in \Z^{n|m}$ by 
\be \label{e69}^L\Gl =^L\Gl_0 + n\lst \mbox{ and }^R\Gl={^R\Gl_0} + m\gl'.\ee
If   $\Gl_0$ and $\Gl$ are as in \eqref{e7} and \eqref{ezr} respectively, it follows that 
\be \label{e67a}
a_i=m(n-i)+ n\gl_{n+1-i}\ee
 and
\be \label{e67b}b_j=n(j-1) + m\gl'_j.\ee
Thus
\be \label{e67d}
a_i - 
a_k=m(k-i)+ n(\gl_{n+1-i}-\gl_{n+1-k})
\ee
and 
\be \label{e67c}b_{j}-b_{k} =n(j-k) + m(\gl'_{j}-\gl'_{k}).\ee
Also
\be \label{e93}
|^L\Gl| - |^L\Gl_0| = n|\gl|, \quad |^R\Gl| -|^R\Gl_0| = m|\gl|.\nn\ee 

\bl \label{abc} If $\Gl=x(\gl)$ is as in \eqref{ezr} and $\ga = \gep_i-\gd_j$ is a positive root, 
\bi
\itema If $\gl\in X_\ga$ then $\Gl\in \Pi_\ga$.
\itemb If $\gl\in X_{-\ga}$ then $\Gl\in \Pi_{-\ga}$.
\itemc If 
$\gl \in X_{\pm\ga}$, then 
$$
x(t_{\pm\ga}(\gl))=\gt_{\pm\ga}(\Gl).$$
\ei
\el
\bpf If $\gl\in X_{\ga}$, then $\gl_{n+1-i} = j-1$ and $\gl'_{j} = n-i$  
by Lemma \ref{rjx}.  Hence $a_i=b_j$ by Equations \eqref{e67a}  and \eqref{e67b}, that is $x(\gl)\in \Pi_\ga$ by  \eqref{e99}.  The proof of (b) is similar.
Finally, if 
$\gl \in X_{-\ga}$, and 
$\mu = t_{-\ga}(\gl)$, then by \eqref{qef}
 \[\mu_k= \left\{ \begin{array}
  {rcl}
  \gl_k -1  & \mbox{if} & k = {n+1-i} \\
\gl_k & \mbox{otherwise} &
\end{array} \right. \]
Similarly
 \[\mu'_k= \left\{ \begin{array}
  {rcl}
  \gl'_k -1  & \mbox{if} & k = j \\
\gl'_k & \mbox{otherwise} &
\end{array} \right. \] It follows from \eqref{e67a}  and \eqref{e67b}
that $x(\mu) = \gt_{-\ga}(\Gl).$ The proof of the other statement is similar.
\epf

\bc \label{aed} 
\bi \itema
If
$\gl \subseteq$ {\bf R}, then $x(\gl)\in   
\cO$. 
\itemb With the notation of \eqref{e7a}, 
$\Gl_0^{(kn)}\in \cO$ for all $k\in \Z$.  
\ei
\ec
\bpf Statement (a) is immediate. For (b), if  $\gl=(m^n)$ then by \eqref{e67a}
 and
\eqref{e67b}, $x(\gl)= \Gl_0  +\bf{mn}$, and this gives the result if $k=1$.    Now take any sequence of morphisms in $\mathfrak  T_{iso}$ that move 
$\Gl_0\in \cO$ to 
$\Gl_0^{(n)}\in \cO$ and apply their inverses in the reverse order to $\Gl_0$ reach 
$\Gl_0^{(-n)}$.  For $|k|>1$ the result follows by induction.  
\epf

\subsection{The coprime case} \label{cpc}
\subsubsection{ Converse to  Lemma \ref{abc} and exceptional cases} \label{spc}
We assume $m,n$ are coprime and obtain partial converses to the statements in Lemma \ref{abc}.  
 Suppose $\ga=\gep_i-\gd_j.$ 
\bl \label{c1} If 
$\Gl=x(\gl)\in \Pi_\ga$, then $\gl\in X_\ga$ unless 
${ \boldsymbol \gl}\subseteq \gl $ and $\ga=\gep_n-\gd_1$. 
\el
\bpf If $\Gl\in \Pi_\ga$, then 
$a_i=b_j$.  
Thus by 
\eqref{e67a}  and \eqref{e67b}
\be \label{gzr} n(\gl_{n+1-i}+1-j)= m(\gl'_j+i-n).\ee 
Since $i\in [m]$, $j\in [n]$, $0\le \gl'_j $ and $\gl_{n+1-i} \le m$, we have  
$$  \gl_{n+1-i}+1-j\le \gl_{n+1-i} \le m \mbox{ and } -n\le \gl'_j-n <  \gl'_j+i-n.$$
Thus, because $m,n$ are coprime  either both sides of \eqref{gzr} are zero, or both equal $mn$.   In the former case  we have $\gl_{n+1-i}=j-1$ and $\gl'_j=n-i$ and 
$\gl\in X_\ga$ by Lemma \ref{rjx}. 
In the latter case  $i=n,$ $j=1$, $\gl_1=m$ and $\gl'_1=n$. 
This means ${ \boldsymbol \gl}\subseteq \gl$ and $\ga=\gep_n-\gd_1$. 
\epf
\bl \label{c2} If  $\Gl\in \Pi_{-\ga}$, then $\gl\in X_{-\ga}$ unless 
$\ga = \gep_1-\gd_{m}$ and $\gl_{n}= \gl'_m=0.$ 
\el
\bpf  
If $\Gl \in \Pi_{-\ga}$ then $a_i-b_j = (n-m)$, which 
by Equations \eqref{e67a}  and \eqref{e67b}
is equivalent to \be \label{gsr} n(\gl_{n+1-i}-j)= m(\gl'_j+i-n-1).\ee 
Both sides of \eqref{gsr} equal 0 or $-mn$.  In the first case 
we have 
$\gl_{n+1-i} = j$  and $\gl'_{j} =  {n+1-i},$ so  
 $\gl\in X_{-\ga}$  by Lemma \ref{rjs}.
In the second,  we deduce $i=1, j=m$ and $\gl_{n}= \gl'_m=0.$\epf
\subsubsection{ Comparison of $\A(\Gl)$ and $\A(\Gl^+)$} \label{spd}
\noi Until the end of Subsection \ref{cpc}, we assume
\be \label{gjd}
\ga = \gep_i-\gd_{j}, \quad \gl\in X_{\ga}, \quad \gl^+=  t_{\ga}(\gl) \quad\mbox{ and } \quad\Gl^+= \gt_{\ga}(\Gl).\ee  
 Then $a_i-b_j =0$, see \eqref{ezr} and \eqref{exy}.  
By Lemma \ref{rjx}, we have\be \label{gjr}\gl_{n+1-i} = j-1 \quad\mbox{ and } \quad\gl'_{j} = n-i.\ee
 Consider the roots $\ga^\ua, \ga^\ur$ which label the boxes immediately above and to the right of box $ \ga$.  These are given by 
$$\ga^\ua =\gep_{i-1}-\gd_{j}, \quad  \ga^\ur=\gep_i-\gd_{j+1}.$$
Of course $\ga^\ua, \ga^\ur$ are only defined if $i>1$ and $j<m$ respectively.  These conditions  will be assumed whenever the notation is used. 
Note that $\Gl^+=\Gl + n\gep_i -m \gd_j$ by \eqref{e9}. 

\bl \label{obc} If $(\Gl^+,\gb) =0$, for $\gb = \gep_i-\gd_k$, then one of the following holds, 
\bi
\itema  $k=j+1$, $\gl'_{j}=\gl'_{j+1}$ and $\gb=\ga^\ur $.
\itemb $\ga = \gep_n-\gd_m$, $\gb = \gep_n-\gd_1$, $\gl \in X^{\red}$ and ${ \boldsymbol \gl}\subseteq \gl^+$. 
\ei
\el
\bpf  Since 
$(\Gl,\ga)=0,$ we have $(\Gl,\gb) =b_{j}-b_{k} =n(j-k) + m(\gl'_{j}-\gl'_{k})$ by \eqref{e67c}. Thus $(\Gl^+,\gb) =0$ iff
\be \label{hjd}
n(j-k+1) = m(\gl'_{k}-\gl'_{j}).\ee
First we show that $k \le j+1$.  Otherwise \eqref{hjd} implies that $k\ge m+j+1$, which is impossible since $j,k\in [m].$ 
If $k=j+1$, then \eqref{hjd} is zero and  we obtain the conclusion in (a). Finally if  $k \le j$, then $j+1=k+m$, which implies that $j=m, k=1, \gl'_{1}=n$  and $\gl'_{m}=0$. Since 
$\gl\in X_{\ga}$, we deduce that $i=n$.
Now
$\gl'_{m}=0$ implies 
$\gl \in X^{\red}$ and 
since 
$\gl'_{1}=n$,
and $\ga = \gep_n-\gd_m$, $\gl^+ = t_{\ga}(\gl)$ contains 
${ \boldsymbol \gl}$.
\epf

\bl \label{pbc} If $(\Gl^+,\gb) =0$, for $\gb = \gep_k-\gd_j$, then one of the following holds, 
\bi
\itema  $k=i-1$, $\gl_{n+1-i}= \gl_{n+2-i}$ 
and $\gb=\ga^\ua $.
\itemb $\ga = \gep_1-\gd_1$, $\gb = \gep_n-\gd_1$, $\gl_{1}=m$, $\gl_{n}=0$ 
 and ${ \boldsymbol \gl}\subseteq \gl^+$. 
\ei
\el
\bpf  As in the proof of the previous Lemma, we have 
using \eqref{e67d},
$(\Gl^+,\gb) =0$ iff
\be \label{hjk}
m(k+1-i) = n(\gl_{n+1-k}-\gl_{n+1-i}).\ee
If $k+1<i $ we obtain 
$i= k+1+n$, which is impossible since $i,k\in [n].$ 
If $i=k+1$, then \eqref{hjk} is zero and  we obtain the conclusion in (a). Finally if  $k+1>i $, then $k+1=i+n$, which implies that $i=1, k=n, \gl_{1}=m$  and $\gl_{n}=0$. Since 
$\gl\in X_{\ga}$, we deduce that $j=1$. The conclusion in (b) follows. 
\epf
\subsubsection{A left inverse to $x$
}\label{daa}  
Let $Z$ be a set of zero entries in $\A(\Gl)$. We say a  descending path $p$ {\it supports } the entries in $Z$ if 
\bi\itema All the entries in $Z$ are above, to the right and adjacent to  $p$. 
\itemb Any other path with the property in (a) lies below and to the left of $p$.
\ei
There is always a unique path which supports more zero entries in $Z$ than any other.  
The entries under the path determine a partition $a(\Gl)$.  The corresponding set of zeroes is denoted  $Z(\Gl)$.  
\\ \\
Assume \eqref{gjd} holds and set $Z(\Gl)^\ga=Z(\Gl)\bsk \{\ga\}$.

\bl \label{le1} The set $Z(\Gl^+)$ satisfies \be \label{fjd}
Z(\Gl)^\ga \subseteq  Z(\Gl^+)\subseteq Z(\Gl)^\ga\cup \{\ga^\ua, \ga^\ur\}.\ee  
and is determined by \eqref{fjd} and the conditions
\bi \itema $\ga^\ur \in Z(\Gl^+)$ iff $\gl'_{j+1}= n-i.$  \itemb $\ga^\ua \in Z(\Gl^+)$ iff $\gl_{n+2-i}= j-1$. 
\ei 
\el

\bpf  
If $k\neq i$ and $\ell\neq j$, then $\A(\Gl)_{k,\ell}= \A(\Gl^+)_{k,\ell}$.  Also since $\A(\Gl)_{i,j} =0$, no other entry in  row $i$ or column $j$ of $\A(\Gl)$ can be zero.  This implies 
$Z(\Gl)^\ga \subseteq  Z(\Gl^+).$  
It remains to consider roots $\gb\neq \ga$ in the same row or column as $\ga$ such that 
$(\Gl^+,\gb)=0.$ We use \eqref{gjr}. If $\gb = \gep_i-\gd_k$, then by
Lemma \ref{obc}, either $k=j+1$, $n-i=\gl'_{j}=\gl'_{j+1}$ and $\gb=\ga^\ur $
or $\gb = \gep_n-\gd_1$. 
However in the latter case, the zero of $\A(\Gl^+)$ in box $\gb$ lies under the path defining 
$\gl^+$, so $\gb \notin Z(\Gl^+)$.  
\\ \\
If $\gb = \gep_k-\gd_j$, then by
Lemma \ref{pbc} either $j-1=\gl_{n+1-i}= \gl_{n+2-i}$ 
and $\gb=\ga^\ua $ or 
$\ga = \gep_1-\gd_1$, $\gb = \gep_n-\gd_1$,
$\gl_{1}=m$, $\gl_{n}=0$ 
 and ${ \boldsymbol \gl}\subseteq \gl^+$. In the latter case, we  can show as above that $\gb \notin Z(\Gl^+)$.  
\epf

\noi For a Young diagram $\gl$, set 
$z(\gl)= \{\gb\in \Gd^+| \gl\in X_\gb\}$ and $z(\Gl)^\ga=z(\Gl)\bsk \{\ga\}$.

\bl \label{le7} The set $z(\gl^+)$ satisfies \be \label{jjd}
z(\gl)^\ga \subseteq  z(\gl^+)\subseteq z(\gl)^\ga\cup \{\ga^\ua, \ga^\ur\}.\ee  
and is determined by \eqref{jjd} and the conditions
\bi \itema $\ga^\ur \in z(\gl^+)$ iff $\gl'_{j+1}= n-i.$  \itemb $\ga^\ua \in z(\gl^+)$ iff $\gl_{n+2-i}= j-1$. 
\ei 
\el
\bpf
The upper boundary of a Young diagram in {\bf R} determines a path.  When $\ga$ is an outer corner  of $\gl$ we consider four consecutive steps in the paths for $\gl$  and $\gl^+$. For example if $\gl'_{j+1}= n-i$ and $\gl_{n+2-i}= j-1$
 we have the paths below with outer corners as shown. 
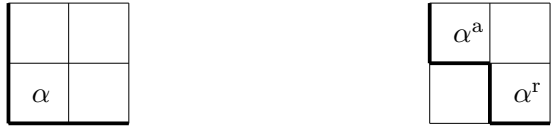
\Bc
\setlength{\unitlength}{0.8cm}
\begin{picture}(10,4)(-1,-.50)
\thinlines
  \linethickness{0.05mm}
\put(9.39,0.35){$\ga^\ur$}
\put(8.39,1.35){$\ga^\ua$}
\multiput(8.0,0)(1,0){3}{\line(0,1){2.0}}
 \multiput(8,0)(0,1){3}{\line(1,0){2.0}}
  \linethickness{0.5mm}
\put(1.0,0){\line(1,0){2}}
\put(1.0,2){\line(0,-1){2}}
  \linethickness{0.05mm}
\put(1.39,0.35){$\ga$}
\multiput(1.0,0)(1,0){3}{\line(0,1){2.0}}
 \multiput(1,0)(0,1){3}{\line(1,0){2.0}}
  \linethickness{0.5mm}
\put(8.0,1){\line(1,0){1}}
\put(9.0,0){\line(1,0){1}}
\put(8.0,2){\line(0,-1){1}}
\put(9.0,1){\line(0,-1){1}}
\end{picture}
\Ec
\noi 
The portion of the path on the left (resp. right) will be denoted d \underline{d} \underline{r}  r (resp. d  \underline{r}   \underline{d} r).  The underlined terms  are edges of $ \fB(\ga)$. 
Passing from the path for $\gl$ to the path in $\gl^+$, the order of the underlined terms is reversed. 
There are four possibilities and they determine $z(\gl^+)$ as in the table below.

\[ \begin{tabular}{|c|c|c|c|c|} \hline
Path in $\gl$
 & r \underline{d} \underline{r}   d& r  \underline{d} \underline{r} r  & d \underline{d} \underline{r} d & d \underline{d} \underline{r} r \\ 
\hline Path in $\gl^+$
& r   \underline{r}   \underline{d} d & r \underline{r} \underline{d} r& d \underline{r}   \underline{d} d& d \underline{r}   \underline{d} r\\
  \hline $z(\gl^+)\bsk z(\gl)^\ga$
&$\emptyset$
&$\{\ga^\ur\}$& $\{\ga^\ua\}$& $\{\ga^\ur, \ga^\ua\}$\\
    \hline
\end{tabular}\]
The result follows from this.
\epf
\noi 
\bl \label{Th1} If  $a(\Gl)=\gl,$ then $\gl^+=a(\Gl^+)$. 

\el
\bpf  
The definition of 
$\gl=a(\Gl)$ depends on the set $Z(\Gl)$ which satisfies the recurrence in Lemma \ref{le1}.
On the other hand the set $z(\Gl)$ from  Lemma \ref{le7} satisfies the same recurrence.
\epf

\bc \label{ch1} If $\gl \subseteq$ {\bf R} and 
$\Gl=x(\gl)$, then $a(\Gl)=\gl.$  Thus $ax$ is the identity map on $X$.
\ec
\bpf  
Obviously the result holds if $|\gl|=0$.  So the result follows by induction on $|\gl|$ 
and Lemma \ref{Th1}. \epf

\subsection{Rotation} \label{rtn}
Suppose $\gl\in X(\da)$.  
If $a_i(\gl)$, $b_j(\gl)$ are the expressions in \eqref{e67a} and 
\eqref{e67b}, then $a_i(\overline{\gl})$, $b_j(\overline{\gl})$ are given by
\be \label{e67f}
a_i(\overline{\gl})=
m(n-i)+ n\overline{\gl}_{n+1-i} = a_{\nu(i)}(\gl)-m \mbox{ for } i\in [n], i\neq 1.\ee
and 
\be \label{e67g}b_j(\overline{\gl})=b_j(\gl) - m \mbox{ for } j\in [m].\nn\ee 
If $i=1$, then since $\overline{\gl}_{\nu(1)} = \overline{\gl}_n =0,$ \eqref{e67f} holds in this case also.  
Hence
\be \label{e67h}
x(\overline{\gl})=
\sum_{i=1}^n  a_i(\overline{\gl})\gep_i - \sum_{j=1}^m 
b_j(\overline{\gl})\gd_j=-{\bf m} +\sum_{i=1}^n  
a_{\nu(i)}(\gl)\gep_i- \sum_{j=1}^m  b_j(\gl)\gd_j.\ee
Now \eqref{e67h} can be written in the form

\be \label{e67k}
x(\overline{\gl})+{\bf m} =\nu(x({\gl})).\ee
The relation \eqref{e67k} can be written in the form 
$\overline{\Gl} 
+{\bf m} =\nu({\Gl}) $ where $
\overline{\Gl} =x(\overline{\gl})$ and $\Gl  =x({\gl}).$  
Thus if $\gl_1=m$ we have an equality of singly augmented matrices

\[\begin{tabular}{|c||ccc|}\hline &     && \\ $^L(\overline{\Gl}
+{\bf m})$
&&$\A(\overline{\Gl})$
 & \\
 & &                              & \\
  \hline
\end{tabular} \quad =\quad 
\begin{tabular}{|c||ccc|}\hline &     &                          & \\ $^L\nu(\Gl)$
&&$\A(\nu(\Gl))$
 & \\
 & &                              & \\
  \hline
\end{tabular}\]
The matrix on the right can be obtained by taking the bottom row of the following matrix
and moving  it to the top, while at the same time moving all other rows to the row below.
\[\begin{tabular}{|c||ccc|}\hline &     &                          & \\ $^L\Gl$
&&$\A(\Gl)$
 & \\
 & &                              & \\
  \hline
\end{tabular}\]
The above matrices have entries in the rectangle obtained from 
{\bf R} by adding  another column on  the left. If  we identify the top and bottom edges of  the expanded rectangle 
 we  obtain a cylinder with circumference $n$ and length $m+1$.  The operation just described using augmented matrices, corresponds to rotating the cylinder through the angle $2\pi/n$. 
Unrolling the text of the cylinder we obtain the graph from Theorem \ref{iir}, see also Example \ref{E3}.
\\ \\ 
 Suppose ${ \boldsymbol \gl}\subseteq \gl =a(\Gl)$. 
By  Lemma \ref{c1} there is no descending path supporting all the zeroes in $\A(\Gl)$.  
However there is a unique descending path supporting all the zeroes in $\A(\nu(\Gl))$.  The partition deternined by this path is 
 $\overline \gl .$  See Example \ref{E2}.

\subsection{The map $x:[X\ti \Z] \lra\Z^{n|m}$.} \label{tmx}
We extend the map $x:X \lra \Z^{n|m}$ from Subsection \ref{coL}.  
Define a map 
$x:X\ti \Z \lra\Z^{n|m}$ by 
\be \label{xdef}
x(\gl,k) = 
\nu^{-k}x({\gl}) +k {\bf m}.\ee
By \eqref{e67k}, 
$x(\overline{\gl},k+1)=x(\gl,k)$ if $\gl\in X(\da)$. 
Thus $x$ extends a map 
$x:[X\ti \Z] \lra\Z^{n|m}$.
\\ \\
If $\Gl= x({\gl},k)$, then 
\be \label{e23}
|^L\Gl| - |^L\Gl_0| = n\deg({\gl},k), \quad |^R\Gl| -|^R\Gl_0| = m\deg({\gl},k).\ee 
Now assume 
\be \label{gjh}
\ga = \gep_i-\gd_{j}, \quad \gb=\nu^{-k}(\ga), \quad\gl\in X_{\pm\ga}, \quad \gl^\pm=  t_{\pm\ga}(\gl) \quad\mbox{ and } \quad\Gl^\pm= \gt_{\pm\ga}(\Gl).\nn\ee

\bl \label{rjg} 
 If $\gl\in X_{\pm\ga},$ and $\Gl= x(\gl)$ then $x(\gl,k)\in 
\Pi_{\pm\gb}$  and 
\be \label{gjt} 
x(\gl^\pm,k)= \gt_{\pm\gb}x(\gl,k).\ee
\el
\bpf If $\gl\in X_{\pm\ga},$ then $\Gl= x(\gl)\in \Pi_{\pm\ga}$ by Lemma \ref{abc}, so  $\nu^{-k}(\Gl)\in  \Pi_{\pm\gb}$ 
and $x(\gl,k)\in \Pi_{\pm\gb}$ by \eqref{xdef}.
By 
\eqref{e9} 
\by \label{e41}  \gt_{\pm\gb}[\nu^{-k}(\Gl) +k{\bf m} ]   
&=&  \nu^{-k}(\Gl) \pm v_{\gb}
+k{\bf m}    \nn\\&=&
\nu^{-k}[\Gl \pm v_{\ga}+k{\bf m} ]
\nn\\&=&
\nu^{-k}[\gt_{\pm\ga}(\Gl)+k{\bf m} ].
\ey  
On the other hand 
\by 
x(\gl^\pm,k)
&=&  \nu^{-k}[x(t_{\pm\ga}(\gl))+k{\bf m} ].
\nn
\ey  
and Lemma \ref{abc} shows that this is equal to \eqref{e41}.  
\epf 

\bt \label{rjt}
The  map 
$x:[X\ti \Z] \lra\Z^{n|m}$ from \eqref{xdef} is injective with image  
$\cO$. 
\et \bpf Any equivalence class has a unique representative in $X^{\red}\ti \Z$.  To show  $x$ is injective suppose $\gl,\mu\in X^{\red}$ and 
\be x(\gl,k) = \nu^{-k}x({\gl}) +k {\bf m}=x(\mu,r) = \nu^{-r}x({\mu}) +r {\bf m}.\ee
We proceed in a sequence of steps.
\bi 
\itema If $k=r$, the result follows from  Corollary \ref{ch1}.  Assume $k>r$ and $d = k-r$.
\itemb 
By Lemma \ref{l73} $\A(\nu^{-k}x({\gl})) = \A(\nu^{-r}x({\mu}))$ and so $\A(x({\gl})) = \A(\nu^{d}x({\mu}))$.  
\itemc  
Also by \eqref{e23}  $|\mu|-|\gl|=dm.$
\itemd 
By repeated use of Lemma \ref{rjg} we can assume 
$r=0,\gl =\emptyset$.  
\iteme 
Now the  unique solution to 
$\Gl_0 =
\A(\nu^{d}x({\mu}))$  is  $\mu = (m^{n-d})$.  But for ${\mu}$ to be in $X^{\red}$ this would imply $d=n.$
\ei
Next we show $x({\gl},k) \in \cO$ for all 
$\gl \subseteq$ {\bf R} and $k\in\Z$.
By Corollary \ref{aed} 
and \eqref{xdef} 
\be \label{ydef}
x(\gl,kn) = 
x({\gl}) +kn {\bf m}\in \cO.\nn\ee
Suppose we have shown  $x({\gl},k-1) \in \cO$ for all 
$\gl \subseteq$ {\bf R}. Then if 
$\gth \subseteq$ {\bf R}  and $\gth_n=0,$ there is a (unique) $\gl \subseteq$ {\bf R} 
such that $\gl_1=m$ and $\overline{\gl} = \gth$.  Then $x({\gth},k) =
x({\gl},k-1) \in \cO$ for all such $\gth$.  Then from \eqref{gjt} it follows that 
$x({\mu},k) \in \cO$
for all
$\mu \subseteq$ {\bf R}. 
\\ \\
On the other hand if $\Gl = x(\gl,k) \in \cO \cap \Pi_{\gb}$, then 
$\gt_\gb(\Gl) \in \cO$ by Lemma \ref{rjg}.  This shows
$\cO\subseteq \Im x$.
\epf \noi 
Since $SV(\gr_{\pm\gb}) =\gt_{\pm\gb}$, we can use \eqref{syx} to write \eqref{gjt} in the form 
\be \label{gqt} 
SV(\gr_{\pm\gb})x[\gl,k] =\gt_{\pm\gb}x[\gl,k]
= x(t_{\pm \nu^k(\gb)}[\gl),k]=x(F(\gr_{\pm\gb})[\gl,k].\nn
\ee
Theorem \ref{iir} follows from this.

\bexa \label{E2}
{\rm 
\noi Take $(n,m) =(2,3)$.  We only indicate the zero entries in $\A(\Gl)$. 
 We start with the  matrix $\hat \A(\Gl_0)$  shown on the left. On the right is the matrix $\hat \A(\Gl_1)$  arising from the morphism corresponding to the unique zero entry in 
$\A(\Gl_0)$. 
\Bc
\setlength{\unitlength}{1cm}
\begin{picture}(10,4)(-1,-.50)
\thinlines
  \linethickness{0.05mm}
\put(1.45,0.35){$0$}
\put(1.45,2.35){$0$}
\put(2.45,2.35){$2$}
\put(3.45,2.35){$4$}
\put(0.45,.35){$0$}
\put(0.45,1.35){$3$}
\multiput(1,0)(-1,0){2}{\line(0,1){3.1}}
\multiput(1.1,0)(1,0){4}{\line(0,1){3.1}}
 \multiput(0,2.1)(0,1){2}{\line(1,0){4.1}}
 \multiput(0,0)(0,1){3}{\line(1,0){4.1}}
  \linethickness{0.5mm}
\put(1.1,0){\line(1,0){3}}
\put(1.1,2){\line(0,-1){2}}
\thinlines
  \linethickness{0.05mm}
\put(9.45,.35){$0$}
\put(8.45,1.35){$0$}
\put(8.45,2.35){$3$}
\put(9.45,2.35){$2$}
\put(10.45,2.35){$4$}
\put(7.45,.35){$2$}
\put(7.45,1.35){$3$}
\multiput(8,0)(-1,0){2}{\line(0,1){3.1}}
\multiput(8.1,0)(1,0){4}{\line(0,1){3.1}}
 \multiput(7,2.1)(0,1){2}{\line(1,0){4.1}}
 \multiput(7,0)(0,1){3}{\line(1,0){4.1}}
  \linethickness{0.5mm}
  \put(8.1,1){\line(1,0){1}}
\put(9.1,1){\line(0,-11){1}}
 \put(9.1,0){\line(1,0){2}}
\put(8.1,2){\line(0,-11){1}}

\end{picture}
\Ec 
Given $\gl \subseteq$ {\bf R} as in Subsection \ref{coL} 
we can use  \eqref{e69}  to find $\Gl=x(\gl) \in \cO$ as in Corollary \ref{aed}.  If 
$\gl=(3,1)$, the vectors $\gl', \lst$ from \eqref{cd} are given by $\gl'=(2,1,1)$ and $ \lst=(1,3)$.  Therefore $\Gl =(5,6|6,5,7).$ 
The matrix $\hat \A(\Gl)$
is given 
on the left below.  There is no descending path supporting all the zeroes, compare Lemma \ref{c1}. The descending path supporting the upper zero
corresponds to the partition $\gl$.  We have $\overline{\gl}=(1,0)$.  
After rotating the diagram on the left, we obtain the diagram on the right where the path corresponds to $\overline{\gl}$.  Note that all entries in this last matrix are the same as those in $\hat \A(\Gl_1)$, except that   the entries in $^L\Gl_1$ and $^R\Gl_1$ have been increased by $m=3$.  In the notation of 
 \eqref{e67k}
$$
x(1,0)+{\bf 3} =\nu(x(3,1)).$$

\Bc
\setlength{\unitlength}{1cm}
\begin{picture}(10,4)(-1,-.50)
\thinlines
  \linethickness{0.05mm}
\put(1.45,0.35){$0$}
\put(2.45,1.35){$0$}
\put(1.45,2.35){$6$}
\put(2.45,2.35){$5$}
\put(3.45,2.35){$7$}
\put(0.45,.35){$6$}
\put(0.45,1.35){$5$}
\multiput(1,0)(-1,0){2}{\line(0,1){3.1}}
\multiput(1.1,0)(1,0){4}{\line(0,1){3.1}}
 \multiput(0,2.1)(0,1){2}{\line(1,0){4.1}}
 \multiput(0,0)(0,1){3}{\line(1,0){4.1}}
\thinlines

  \linethickness{0.05mm}
\put(8.45,1.35){$0$}
\put(9.45,.35){$0$}
\put(8.45,2.35){$6$}
\put(9.45,2.35){$5$}
\put(10.45,2.35){$7$}
\put(7.45,.35){$5$}
\put(7.45,1.35){$6$}
\multiput(8,0)(-1,0){2}{\line(0,1){3.1}}
\multiput(8.1,0)(1,0){4}{\line(0,1){3.1}}
 \multiput(7,2.1)(0,1){2}{\line(1,0){4.1}}
 \multiput(7,0)(0,1){3}{\line(1,0){4.1}}
  \linethickness{0.5mm}

  \put(8.1,1){\line(1,0){1}}
\put(9.1,1){\line(0,-11){1}}
 \put(9.1,0){\line(1,0){2}}
\put(8.1,2){\line(0,-11){1}}
  \put(1.1,2){\line(1,0){1}}
\put(2.1,2){\line(0,-11){1}}
 \put(2.1,1){\line(1,0){2}}
\put(4.1,1){\line(0,-11){1}}

\end{picture}
\Ec
}\eexa

\noi

\section{Borel subalgebras} \label{akmls} 

\subsection{Shuffles and Borel subalgebras of  $\fgo=\fsl(n|m)$} \label{fte}
In this Section $\fgo$ will denote the Lie superalgebra $\fsl(n|m)$ with $m>n$ or $\fgl(n|n)$.   Set $r=m+n-1$  in the former case and $r=2n-1$ in the latter. 
Let $\fho$ be the subalgebra of diagonal matrices in  $\fgo$. 
Suppose that $ A^0  = (a^0_{ij})$
is the Cartan matrix of $\fgo$ corresponding to the shuffle $\gs$. 

\subsubsection{Outline} \label{spafn}   
For simplicity we assume here that $\fgo=\fsl(n|m).$  
Let $\pigo$ be the distinguished set of simple roots for $\fgo$, $\fbo$ the distinguished Borel subalgebra   
and let  $\gb=\gep_n-\gd_1$ be the unique odd root in $\pigo$.  For any shuffle $\gs$ we have a corresponding Borel $\fbo(\gs)$ as in \eqref{bji} with simple roots 
$\pigo(\gs) $
and  highest root $\ghs$.
Set $\ga_0(\gs)=\ogd -\ghs$.  In the affinization $\Lgh$ set 
$ \Pi(\gs) =\pigo(\gs)\cup\{\ga_0(\gs)\}.$  The corresponding extended Dynkin-Kac diagram is obtained from the Dynkin-Kac diagram for $\fgo$
arising from the simple roots $\pigo(\gs) $ 
by adjoining a node corresponding to the root $\ga_0(\gs)$.   If 
$\gs(m')= n$, then  
by deleting the node of the extended diagram 
corresponding to the simple root $\gep_{\gs(r)}-\gep_{\gs(r+1)}$ we obtain a set of simple roots for a
different subalgebra $\fgo(1)$ of $\Lgh$ which is isomorphic to $\fgo$.  This allows us to perform further odd reflections that do not arise from $\fgo$.    
Below we provide further details of this procedure. In particular we need generators for the algebras to explain how $\fgo(1)$ is embedded in  $\Lgh$. 
\subsection{$\fgo$ as a contragredient Lie superalgebra} \label{afn1}
\noindent The algebra  $\fgo$ can be constructed as a contragredient Lie superalgebra   $\fgo = \fg(A^0, \tau^0)$. 
We briefly recall some details. Let  ${\stackrel{{\rm o}}{\fb}}{}(\gs)$ be the Borel subalgebra of $\fgo$ corresponding to the shuffle $\gs$. We are given subsets 
\be \label{pihid}\pigo = \{\alpha_1, \ldots, \alpha_r\} \;\; \mbox{and} \;\;
\stackrel{\rm o}{\Pi^\vee} = \{H_1,\ldots, H_r\}.\ee
of $\fho$ and $\fho^*$ respectively.  The set 
$\pigo$ is the set of simple roots of ${\stackrel{{\rm o}}{\fb}}{}(\gs)$  and the $\ga_i$ are given by $\ga_i=\gep_{\gs(i)}-\gep_{\gs(i+1)}$. For $i\in [r]$ let   $E_i, F_i$ be root vectors with weights $\ga_i,-\ga_i$  respectively and set    $H_i=[E_i, F_i]$.  Then $\fgo$ is generated by the $E_i, F_i$ and the Cartan subalgebra $\fho$ of diagonal matrices. These elements satisfy the relations
(18.2.2)-(18.2.4) from \cite{M101}.
The subalgebra $\fb(\gs)$ is generated by the $E_i $ and  $\fho.$  
The  matrix $A^0$ has entries  $a^0_{ij}=\alpha_i(H_j)$.  
Also all elements of $\fho$ are even, and for $ i \in [r]$, $E_i$ and $F_i$ are odd if and only if $i \in \tau^0$.

\bl \label{pih} A minimal realization of  $A^0$ is given by   
 $(\fho, \pigo, \pigov)$.  We have $\fgo \cong \fg(A^0, \tau^0)$.
\el \bpf If 
$\fgo =\fsl(n|m)$, this is contained in \cite{M101} Theorem 5.3.5. The result for 
$\fgl(n|n)$ is also well-known.  For the distinguished set of simple roots  it is \cite{M101}, Exercise 5.6.12.
\epf
\noi 
If we perform an odd reflection using an odd simple root in $\pigo$ we obtain a new set of generators for $\fgo$ and a new Cartan matrix $A_1^0$.  The relationship between $A^0$ and $A^0_1$ is explained in \cite{M101} Section 3.5 and  Lemma 5.3.4.
\subsection{Affinization} \label{afn}
Starting from $A^0$ and $\theta$ we construct
another matrix $A = (a_{ij})$ with an extra row and column labeled
by 0. 
 There is a unique vector $v$ such that the augmented matrix below $A$ has all row sums and column sums equal to zero.  We refer to the  extra row and column as row and column 0.  
\be \label{mtr} A=\begin{tabular}{|c||ccc|} \hline
0&  &$v^\trans$&
\\ \hline\hline &     &                          & \\ $v$
&&$A^0$
 & \\
 & &                              & \\
  \hline
\end{tabular}\ee
If $\gth$ is even (resp. odd), let $\tau = \tau^0 $ (resp. $\tau = \tau^0 \cup \{ 0 \}$). We give an
 explicit construction of the contragredient Lie superalgebra
 $\fg(A,\tau)$.\\
\\
 The {\it loop algebra} $\Lgo$ is the Lie superalgebra with graded
 components, $(i = 0,1)$
 \[ \Lgo_i = \ttk[t,t^{-1}] \otimes \fgo_i.\]
Next the 
invariant bilinear form 
 $( \; , \;)$ 
on $\fgo$ is used to form a
 one dimensional central extension $\Lgo \oplus \ttk c $ of $\Lgo$, \cite{M101} Lemma 18.2.1. The product is given by 
\be \label{mjr} [t^m \otimes x, t^n \otimes y] = t^{m+n} \otimes [x,y] + m(x,y)
\delta_{m,-n} c. \ee
Define a derivation $d$ on $\Lgo \oplus \ttk c$ such that $d(c)=0$ and $d$ acts as 
$t\frac{d}{dt}$ on the loop algebra. 
The {\it affinization} of $\fgo$ is the Lie superalgebra 
 \[ \Lgh = \Lgo \oplus \ttk c \oplus \ttk d. \]
 We make $\Lgh$ into a Lie superalgebra with $d$ even, such
 that the product on $\Lgh$ extends the product already defined on
$\Lgo \oplus Kc,$ by requiring that
$ [d,x] = d(x)$
for all $x \in \Lgo \oplus Kc$.
Let  $\theta =
 \sum^r_{i=1}\alpha_i$  be the highest root of $\fgo$ and choose elements 
$F_0 \in \; \fgt, E_0 \in \;
\fgm,$ suitably normalized as in \cite{M101} Equation (18.2.12) and such  that $(E_0 ,F_0)=1$.  
Set $$\fh = (1 \otimes \fho) \op \ttk c \oplus \ttk d,$$ and extend
the elements $\alpha_1, \ldots, \alpha_r \in \; \fhs$ to linear
forms on $\fh$ by setting $\alpha_i(c) = \alpha_i(d) = 0$ for $ i \in  [r]$. Define $\overline{\delta} \in \fh^*$ by
\[ \overline{\delta}(c) = 0, \;\; \overline{\delta}(d) = 1, \;\; \overline{\delta}(h) = 0 \;\; \mbox{for} \;\; h
\in \; \fho . \] 
In  $\Lgh$, set 
\[ e_i = 1 \otimes E_i, \quad f_i = 1 \otimes F_i, \quad h_i = 1 \otimes H_i. \]
Set $e_0 = t \otimes E_0, f_0
= t^{-1} \otimes F_0,$ and
\be\label{efo}
 h_0 = [e_0, f_0]  =  [E_0, F_0] + (E_0, F_0)c \nonumber \\
           =  -(1 \otimes h_\theta)+ c.
\ee
We set $\alpha_0 = \overline{\delta} - \theta \in \; \fh^*$,  and $h_i = 1 \otimes H_i \in \fh$ for $i \in [r]$.
 Then set 
\be\label{187}
\Pi =
 \{\alpha_0, \alpha_1, \ldots, \alpha_r\}, \quad\mbox{ and } \quad
 \Pi^\vee =
 \{h_0, h_1, \ldots, h_r\}.\ee  Now suppose 
$\fgo =\fsl(n|m)$. 
 Then the triple $(\fh, \Pi, \Pi^\vee)$ is a minimal realization of $A$ and  $\Lgh \cong \fg(A,\tau)$,
\cite{M101}, Lemma 18.2.2 and we have by  \cite{M101} Theorem 18.2.5. 
  \bt \label{185} With the above notation, $\Lgh \cong \fg(A,\tau)$.  This algebra is generated by $e_0, e_1,
\ldots, e_r,$ $ f_0, f_1, \ldots, f_r$ and $\fh$.
\et \noi 
We explain how to modify Theorem \ref{185} in the case where 
$\fgo =\fgl(n|n).$  In this case $A$ is a $2n\ti 2n$ matrix 
of rank $2n-2$.  Using the fact that $\dim \fh = 2n+2$ it is easy to show that 
$(\fh, \Pi, \Pi^\vee)$ is a minimal realization of $A$.
Let $\tilde \fg=\bop_{i\in \Z} \tilde \fg $ be the graded 
subalgebra of $\Lgh$ given by 
\[\tilde \fg(i) =\left\{ \begin{array}
  {ccl}
 \fgl(n|n) \op \ttk c \op \ttk d& \mbox{if} & i=0\\
\fsl(n|n)t^i& \quad \quad \mbox{otherwise} &
\end{array} \right. \]
 This algebra appears in \cite{vdL} equation 6.1.
\bt \label{105}  The subalgebra of $\Lgh$ generated by $e_0, e_1,
\ldots, e_r, f_0, f_1, \ldots, f_r$ and $\fh$ is   
$\tilde \fg$. 
\et
 \bpf We indicate the changes that need to be made to the proof of \cite{M101} Theorem 18.2.5.
Let $L$ be the subalgebra of $\Lgh$ generated by $e_0, e_1,
\ldots, e_r, f_0, f_1, \ldots, f_r$ and $\fh$.  We claim that $L =
\tilde \fg$. Since $\fgo$ is generated by $\fho$ and $e_1, \ldots, e_r,$ $
f_1, \ldots, f_r$, it follows that $1 \otimes \fgo \; \subseteq L$.
Let
\[ \mathfrak{a} = \{ x \in \; \fgo | t \otimes x \in L \}. \]
Since $E_0 \in \fa,$ 
$\mathfrak{a}$ is an ideal of $\fgo$ not contained in the center. 
The claim follows from this.
\epf

\brs\label{186} {\rm \bi \itema In the terminology of \cite{GHS}, Theorem 9.1.1 $\tilde \fg$ is the universal root algebra of type $A(n-1|n-1)^{(1)}$.
\itemb
The proof of Theorem \ref{185} shows more. The Cartan subalgebra $\fh$ appears in the definition of both $\Lgh$ and $\fg(A,\tau)$ and  the isomorphism is the identity on $\fh$. Also  the generators $e_0, e_1, \ldots, e_r$ and $f_0, f_1, \ldots, f_r$  of $\Lgh$ correspond under the isomorphism to the canonical generators of $\fg(A,\tau)$  as a contragredient Lie superalgebra. If $\fgo=\fgl(n|n)$ it is convenient to make the non-standard definition that $\fg(A,\tau)$ is the subalgebra 
$\tilde \fg$ of $\Lgh$ from Theorem \ref{105} and we refer to this algebra as the {\it affinization} of $\fgo$. Note however that the maximal ideal $\fr$ of $\tilde \fg$ such that $\fr\cap\fh =0$ is non-zero.  It is given by  
$\fr=\bop_{i\in \Z\bsk \{0\}}It^i$ where $I\in \fsl(n|n)$ is the identity matrix.
\ei
}\ers

\subsection{A fundamental construction} \label{afc}
From now on we assume $\fg$ is the affinization of $\fg$.  Thus 
$\fg=\Lgh$ or 
$\tilde \fg$ if  $\fgo=\fsl(n|m)$ or  $\fgl(n|n)$ respectively.  
Let $A^0(1)$ be the submatrix of the matrix $A$ from \eqref{mtr} located in the upper left $r\ti r$ corner.  There is a unique vector $w$ such that $A$ can be written in the form
\[A=\begin{tabular}{|ccc||c|} \hline
&  &&\\
 &    $A^0(1)$  &                          & $w$\\ 
&& & \\
 \hline\hline & $w^\trans$&                              & 0\\
  \hline
\end{tabular}\] 
Let $$\tau^0(1) = \tau \cap \{ 0,\ldots, r-1 \},  \;
{\stackrel{{\rm o}}{\Pi}}(1)= \{\alpha_0, \ldots, \alpha_{r-1}\}, \; 
{\stackrel{{\rm o}}{\Pi}^\vee}(1) =
 \{h_0, \ldots, h_{r-1}\} \; \mbox{and} \; {\stackrel{\rm o}{\fh}}(1) =\span \pigo^\vee(1) .$$ 
Then the triple $({\stackrel{\rm o}{\fh}}(1), 
{\stackrel{{\rm o}}{\Pi}}(1),{\stackrel{\rm o}{\Pi^\vee}}(1))$ is a minimal realization of $A^0(1)$.  If $
\stackrel{\rm o}{\fg}(1) = \fg(A^0(1), \tau^0(1))$, then  by  Theorem \ref{185} we have ${\widehat{L}(\stackrel{\rm _o}{\fg}}(1))\cong \fg(A,\tau)$.  
Moreover viewed as a subalgebra of $\fg(A,\tau)$ by means of this isomorphism, ${\stackrel{{\rm o}}{\fg}}(1)$ is generated by $ \{e_0, \ldots, e_{r-1}\},  \{f_0, \ldots, f_{r-1}\}$ and ${\stackrel{\rm o}{\fh}}(1)$, compare Remark \ref{186}.
Now assume  that $\gl=\gz(\gs)$ has $\gl_1=m$, equivalently $\gs(m')= n$ by  Lemma \ref{rags}. 
Then  ${\stackrel{{\rm o}}{\fg}}(1)$ is  the subalgebra of $\Lgh$ whose 
Dynkin-Kac diagram is obtained from the diagram 
for $\Lgh$ by deleting the node corresponding to the simple root $\gep_{\gs(r)}-\gep_{\gs(r+1)}$ as explained in Subsection \ref{spafn}.  
Then ${\stackrel{{\rm o}}{\fg}}(1)\cong \fgo$. 
\noi We say that the Borel subalgebra of $\Lgh$ with simple roots $\Pi$ as in \eqref{187} {\it is the extension to} $\Lgh$  of the Borel subalgebra of ${\stackrel{{\rm o}}{\fg}}$ with simple roots ${\stackrel{{\rm o}}{\Pi}}$ as in \eqref{pihid}.  Note that the {\it extending root} $\alpha_0 = \overline{\delta} - \theta$ is determined by ${\stackrel{{\rm o}}{\Pi}}.$  In a similar way if  $\gl=\gz(\gs)$ is a partition with $\gl_1 = m$, then $\fb$ is the extension to $\Lgh$  of the Borel subalgebra of ${\stackrel{{\rm o}}{\fg}}(1)$ with simple roots ${\stackrel{{\rm o}}{\Pi}}(1)$.
\\ \\
Next we explain this procedure in terms of shuffles.   Suppose  $\gs\in 
\Sh$ satisfies $\gs(r) = n$. The set of simple roots $
{\stackrel{{\rm o}}{\Pi}}(\gs)$  of is given by \eqref{pgs} and the extending root is $\overline{\gd}+ \gep_{n}- \gd_1$. Therefore 
\be\label{efx}
 \stackrel{{\rm o}}{\Pi}(1) = \{ \overline{\gd}+ \gep_{n}- \gep_{\gs(1)},\gep_{\gs(i)}- \gep_{\gs(i+1)}|i\in [r-1]\}.
\ee
However $\fgo(1) \cong \fsl(n|m)$ or $\fgl(n|n)$ and under this isomorphism, the set of simple roots in \eqref{efx} corresponds to the shuffle $\overline{\gs}$. 
Repeating this process we obtain subalgebras 
${\stackrel{{\rm o}}{\fg}}(k)$ of $\Lgh$ for $k\in  \Z$, all isomorphic to $\fgo$ and Borel subalgebras ${\stackrel{{\rm o}}{\fb}}{}(\gs,k)$ of 
${\stackrel{{\rm o}}{\fg}}(k)$ such that the extensions of ${\stackrel{{\rm o}}{\fb}}{}(\gs,k)$
and ${\stackrel{{\rm o}}{\fb}}{}(\overline{\gs},k+1)$ to  $\Lgh$ are equal, whenever $\gl=\gz(\gs)$ satisfies $\gl_1=m$. For the notation see Lemma \ref{ymn}. 
Let $\cB$ be the set of Borel subalgebras of $\Lgh$ that are extensions of Borel subalgebras of Borels in the subalgebras $\fgo(k)$ for  $k\in \Z.$  
From this discussion, we have
\bl \label{crb} Two pairs $({\stackrel{{\rm o}}{\fb}}{}(\gs),k)$ and $({\stackrel{{\rm o}}{\fb}}{}(\gt),j)$  in ${\stackrel{{\rm o}}{\cB}} \ti \Z$ are equivalent  under $\sim$ iff ${\stackrel{{\rm o}}{\fb}}{}(\gs,k)$  and ${\stackrel{{\rm o}}{\fb}}{}(\gt,j)$ extend to the same Borel subalgebra of $\fg$.  Thus there is bijection ${\stackrel{{\rm o}}{\cB}} \ti \Z \rl \cB$ where $[{\stackrel{{\rm o}}{\fb}}{}(\gs),k]$  corresponds to the Borel subalgebra of $\fg$ extended from ${\stackrel{{\rm o}}{\fb}}{}(\gs,k)$.
\el
\bt \label{srb} The groupoid  $
\mathfrak  T_{iso}$ acts on $\cB$ by odd reflections and we have a commutative 
diagram 
\[
\xymatrix@C=2pc@R=1pc{
\mathfrak  T_{iso}\ar@{<->}_= [dd]&&
\fS(\cB) \ar@{<-}[ll]: \ar@{<->}^\cong[dd]&\\ \\
\mathfrak  T_{iso} \ar@{-}[rr] &&
\fS[{\stackrel{{\rm o}}{\cB}} \ti \Z]\ar@{<-}[ll] &}\]
\et
\bpf It follows from \eqref{kxz} that the correspondence from Lemma \ref{crb} is compatible with odd reflections.  
\epf \noi 
It is important to point out that the permutation $\nu$ is involved in passing from  the shuffle $\gs$ to $\overline{\gs}$, see Lemma \ref{ymn}. This is necessary if we are to view ${\stackrel{{\rm o}}{\fb}}{}(\gs,k)$
and ${\stackrel{{\rm o}}{\fb}}{}(\overline{\gs},k+1)$ as subalgebras of $\fgo(k)$ and $\fgo(k+1)$.  We refer to ${\stackrel{{\rm o}}{\fb}}{}(\gs,k)$ as the {\it local name } for this Borel subalgebra.  Likewise for all the systems of simple roots we have given so far, we have used their local names. However if we want to express the simple roots of ${\stackrel{{\rm o}}{\fb}}{}(\gs,k)$ in terms of the distinguished roots and the extending root  of $\fgo(0)$, we use {\it global names}, which do not involve $\nu$, see Example \ref{E5}.  In the next Subsection, we explain how this can be done in several  important cases.

\subsection{Simple Roots}\label{sro}
If $\fgo =\fgl(n|n)$, take $m=n$ in what follows. We write a set of  simple roots of $\fgo$ as an ordered set.  The first element of the set corresponds to the leftmost node of the Dynkin-Kac diagram.
\subsubsection{Simple  roots of $\fgo$}\label{sro}
The distinguished set of simple roots of $\fgo$ is
\be \label{f5}   \gb_{ i} = \left\{ \begin{array}
  {rcl}\gep_{i}-\gep_{i+ 1}
& \mbox{if} & 1\le i \le n-1 \\
\gep_{n}-\gd_{ 1}& \mbox{if} & i=n\\
\gd_{i-n}-\gd_{i+1-n}
& \mbox{if} & n+1\le i \le n+m-1.
\end{array} \right. \ee
The extending root is 
$\overline{\gd}-\sum_{ i=1}^{m+n-1}\gb_{i} .$
\bl \label{f6}   
The anti-distinguished set of simple roots of $\fgo$ is
\by \label{f9}  && (\gd_{ 1}- \gd_{2}, \;\ldots \; \gd_{m-1}-\gd_{m}, \gd_{m}-\gep_{1}, \;\ldots \;   \gep_{ n-1} - \gep_{ n})\nn\\
&=& (\gb_{n+ 1},  \ldots, \gb_{m+n-1}, - \sum_{ i=1}^{m+n-1}\gb_{i},\gb_{ 1}, \ldots, \gb_{n-1}).\nn\ey 
The extending root is 
$\overline{\gd}+ \gb_{n}.$
\el

\subsubsection{Two procedures}\label{2pr}
We can find the roots of ${\stackrel{{\rm o}}{\fb}}{}(\gs,k)$ by combining two procedures.  The first uses odd reflections and corresponds to adding  an outer corner
to    a Young diagram. For a description of the simple roots in terms of shuffles see \eqref{pgs}. The second involves passing from a Borel in  $\fgo(k)$ to another in $\fgo(k+1)$ such that both extend to the same  Borel in a common affinization.  This 
corresponds to adding or deleting a row in a Young diagram. This procedure has a simple combinatorial interpretation.
Suppose the simple roots of the Borel subalgebra  ${\stackrel{{\rm o}}{\fb}}{}(\gs,k)$ are given as an ordered sequence $ (\ga_1, \ldots, \ga_{r})$ 
where $\ga_i$ is the weight of the root vector $e_i$ in  Subsection \ref{afn}.  The extending root is
 $\ga_r=\overline{\delta} - \sum_{ i=1}^{r}\ga_{i}$.  Then the simple roots of the Borel subalgebra ${\stackrel{{\rm o}}{\fb}}{}(\overline{\gs},k+1)$ are given by 
$(\ga_0, \ldots, \ga_{r-1})$. 

\bl \label{f61}   
Suppose the distinguished set of simple roots of $\fgo(j)$ is
\by \label{f9} (\gb_{1},  \ldots, \gb_{n-1},\gb_{n}, \gb_{n+ 1},  \ldots, \gb_{r} ).\ey 
so that the unique odd simple root is 
$\gb_{n}.$ 
Then for $j\in [n]$, 
the set of roots of $\fgo(j)$ corresponding  to the partition $(m^j,0^{n-j})$ is
\by \label{z91}  (\gb_{1},  \ldots, \gb_{n-j-1},
 \sum_{ i=n-j}^{n}\gb_i, \gb_{n+ 1},  \ldots, \gb_{r},- \sum_{ i=n+1-j}^{r}\gb_i, \gb_{n+1-j}, \ldots, \gb_{n-1}).\ey 
If $j=n$ the terms $\gb_{n+ 1},  \ldots, \gb_{r}$ and $\sum_{ i=n-j}^{n}\gb_i$ are not present in \eqref{z91}. For $j=n$ we obtain the anti-distinguished set of simple roots.
\el
\bpf  Fix $j$ with $0\le j \le n-1$.  To pass from the Borel  corresponding  to $(m^j,0^{n-j})$  to the Borel corresponding  to  the partition $(m^{j+1},0^{n-j-1})$ perform a sequence of odd reflections starting with $\gc =
 \sum_{ i=n-j}^{n}\gb_i$ followed in order by $\gc+ \gb_{n+ 1}, \ldots, \gc + \sum_{ i=n+1}^{\ell}\gb_i, \ldots, \gc + \sum_{ i=n+1}^{r}\gb_i.$  We leave it to the reader to check the details.
\epf\noi 
It is convenient to concatenate ordered sets of roots of the form  

$$\Xi_j = (\gep_{ 1}- \gep_{2}, \;\ldots \; \gep_{ \ell}- \gep_{\ell+1}, \;\ldots \;, \gep_{n-j-1}-\gep_{n-j}),$$
$$\Psi_j = (\gep_{n+1-j}-\gep_{n+2-j}, \ldots 
\gep_{ \ell}- \gep_{\ell+1}, \;\ldots \; ,\gep_{n-1}-\gep_{n})$$ and 
\be\label{g112} 
\nabla=(\gd_{ 1}- \gd_{2}, \;\ldots \; ,\gd_{m-1}-\gd_{m}
).\ee\noi 
Take $k= 0$ in Lemma \ref{f61} and define the roots $\gb_i$ as in \eqref{f5}.
Then we can rewrite \eqref{z91} as
\be\label{g12} 
( \Xi_j,  \sum_{ i=n-j}^{n}\gb_i, \nabla,
- \sum_{ i=n+1-j}^{r}\gb_i,
\Psi_{j}).\ee If $j=n$ (resp. $j=0$) the first two (resp. last two) terms are not present in \eqref{g12}. For $0<j<n$ the extending root is $\overline{\gd}+\gep_n -\gep_{1}$. 
\bl\label{g16} For $0\le j<n$ the distinguished set of roots for $\fgo(j)$ is  
\be\label{g15} 
(\Psi_{j},  \overline{\delta} +\gep_{n} - \gep_{1},\Xi_j,  \gep_{n-j}-\gd_1, \nabla).
\ee 
\el
\bpf  
Apply the second procedure described in Subsection \ref{2pr} $j$ times to the set of roots in \eqref{g12}.
\epf\noi 
\bl\label{g17} If the anti-distinguished set of roots for 
$\fgo(j)$, is
$(\gc_{1},  \ldots, \gc_{m},   \ldots, \gc_{r} )$, with
$\gc_{m}$  the unique odd simple root,
then the  distinguished set of roots for $\fgo(j+n)$ 
is 
$$(\gc_{m+1},  \ldots, \gc_{n},   \ldots, \gc_{r},\overline{\delta} - \sum_{ i=1}^{r}\gc_i,\gc_{1},  \ldots, \gc_{n},   \ldots, \gc_{m-1} ).$$
\el
\bpf Take  $j=n$ in \eqref{z91} or \eqref{g12} and apply the same procedure $n$ times.
\epf
\subsubsection{Distinguished and  anti-distinguished roots of $\fgo(kn)$}\label{Pad}
\bl \label{qxc} 
\bi
\itema
The distinguished set of simple roots of $\fgo(kn)$ is 
$$(\gep_{ 1}- \gep_{2}, \;\ldots \; ,\gep_{n-1}-\gep_{n},  
k\overline{\gd}+ \gep_{n}- \gd_1, \;\ldots \;   \gd_{ m-1} - \gd_{m}).$$
The extending root is 
$-(k-1)\overline{\gd}+\gd_m -\gep_{1}.$
\itemb 
  For $k\in \Z$, the anti-distinguished set of simple roots of $\fgo(kn)$ is
\be\label{g111} 
(\gd_{ 1}- \gd_{2}, \;\ldots \; ,\gd_{m-1}-\gd_{m},  
\gd_{m}-\gep_{1}-k\overline{\gd}, \;\ldots \;   \gep_{ n-1} - \gep_{ n}).\ee
The extending root is 
$(k+1)\overline{\gd}+\gep_{n}- \gd_1.$
\ei\el
\bpf We know (a) holds when $k=0$.  Also statement (a) implies (b) by 
 Lemma \ref{f6} with $k\overline{\gd}+ \gep_{n}- \gd_1$ playing the role of $\gb_n$.  Given (b) we obtain (a) with $k$ replaced by $k+1$ from Lemma \ref{g17}.  This proves (a) and (b) for non-negative $k$. Negative values of $k$ can be handled in a similar way.
\epf
 
\subsubsection{The largest equivalence class}\label{lec}
\bl \label{qey} \bi \itemo
\itema  Every equivalence class in  $X\ti \Z$ contains a unique representative of the form $(\gl,k)$ with ${ \boldsymbol \gl} \subseteq  \gl$ 
and a unique representative of the form $(\mu,\ell)$ with 
$\mu \in X^{\str}$. 
\itemb  Set $d=\ell-k\ge 0$,  $\gl^{(0)}=\gl$ and, if $d>0$ $\gl^{(i)}=\overline{\gl}^{(i-1)}$ for $i\in [d]$.  Then the equivalence class of $(\gl,k)$ is $\{(\gl^{(i)},k+i)| 0\le i\le d\}.$
\ei
\el  \noi
By Lemma \ref{qey} any equivalence class under $\sim$ can contain at most $n+1$ elements.  The  maximum is achieved when
 $\gl=(m^n)$. Any other equivalence class  can contain at most $n$ elements. Set $\gl^{(i)}=(m^{n-i}, 0^{i})$, then  since $\overline{\gl}^{(i-1)}= {\gl^{(i)}}$ 
for $i\in [n]$ we have 
\be\label{h11}  (\gl^{(0)},0) \sim \ldots \sim (\gl^{(i)},i)\sim \ldots  \sim (\gl^{(n)},n).\ee
Write $\gl^{(i)}= \gz(\gs^{(i)})$ as in Lemma \ref{xmn}.  
Then 
\be\label{f11}  (\fbo(\gs^{(0)}),0) \sim \ldots \sim (\fbo(\gs^{(i)}),i)\sim \ldots  \sim (\fbo(\gs^{(n)}),n).\ee
It follows from 
 Equation \eqref{f11} and Lemma \ref{crb} that there is a Borel subalgebra $\fb$  of $\Lgh$ that is the extension of the Borel $\fbo(\gs^{(i)},i) $
in each subalgebra $\fgo(i)$ of $\Lgh$ for $0\le i\le n$.
The simple  roots 
$\Upsilon$ of $\fb$ are obtained by adjoining the extending root 
$\overline{\gd}+ \gb_{n}$ to the set of roots from Lemma \ref{f6}.  
Thus 
\be \label{f19}\Upsilon= (\overline{\gd}+ \gb_{n}, \gd_{ 1}- \gd_{2}, \;\ldots \; \gd_{m-1}-\gd_{m}, \gd_{m}-\gep_{1}, \;\ldots \;,   \gep_{ n-1} - \gep_{ n}).\ee
For $i\in [n-1]$ we obtain the simple roots of $\fbo(\gs^{(i)},i) $ by removing the root $\gep_{ n-i} - \gep_{ n+1-i}$ from the set $\Upsilon$. Removing $\gd_{m}-\gep_{1}$ gives the simple roots of $\fbo(\gs^{(n)},n) $, which is the distintinguished Borel in $\fgo(n)$, and of course removing $\overline{\gd}+ \gb_{n}$ we recover the simple roots from 
Lemma \ref{f6}. 

\bexa \label{E4}{\rm We consider certain Borel subalgebras of $\Lgh$, when $\fgo=\fsl(2|3)$. 
Let $\gs, \gs_1$ be the shuffles with
$${\it{{\underline{\sigma}}}} = ( 1',1, 2',3', 2 ) \mbox{ and } {\it{{\underline{\sigma}}}}_1 = ( 1',2', 1, 3', 2).$$
Let $\fb(\gs)$ and $\fb(\gs_1)$ be the Borel subalgebras of ${\stackrel{{\rm o}}{\fg}}$ corresponding to $\gs$ and $ \gs_1$.  
Note that $\fb(\gs)$ and $\fb(\gs_1)$  are related by odd reflections using the roots
$\pm(\gep_1-\gd_2).$ We give two examples of the procedure defined above.  In the first
(resp. second) the Dynkin-Kac diagram for $\fb(\gs)$ (resp.  $\fb(\gs_1)$) is given on the left, the corresponding diagram for the affinization is in the middle and if $\gt$ (resp. $\gt_1$) is the shuffle obtained by the procedure, the Dynkin-Kac diagram for $\fb(\gt)$ (resp.  $\fb(\gt_1)$) is given on the  right. 
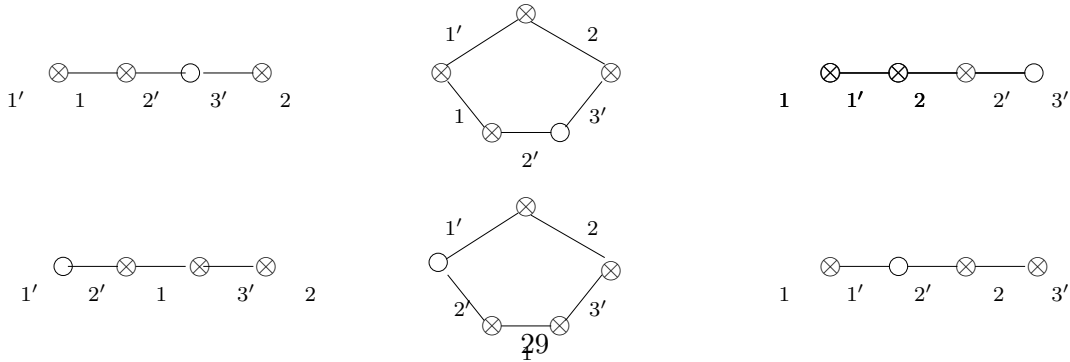
\Bc
\setlength{\unitlength}{0.9cm}
\begin{picture}(10,1)(-0.5,1.70)
\thinlines
  \linethickness{.09mm}
\put(-2.4,1.81){\line(1,0){0.74}}
\put(-1.4,1.81){\line(1,0){0.74}}
\put(-0.4,1.81){\line(1,0){0.74}}

\put(-3.3,1.3){$\scriptstyle 1'$}
\put(-2.3,1.3){$\scriptstyle 1$}
\put(-1.3,1.3){$\scriptstyle 2'$}%
\put(-0.3,1.3){$\scriptstyle 3'$}
\put(.73,1.3){$\scriptstyle 2$}%
\put(0.3,1.7){{$\otimes$}}

\put(-2.7,1.7){{$\otimes$}}
\put(-1.7,1.7){{$\otimes$}}\put(-0.59,1.81){\circle{0.27}}
\put(9,1.81){\line(1,0){0.74}}
\put(10,1.81){\line(1,0){0.74}}
\put(11,1.81){\line(1,0){0.74}}

\put(8.1,1.3){$\scriptstyle 1$}
\put(9.1,1.3){$\scriptstyle 1'$}
\put(10.1,1.3){$\scriptstyle 2$}
\put(11.27,1.3){$\scriptstyle 2'$}\put(12.13,1.3){$\scriptstyle 3'$}%

\put(8.7,1.7){{$\otimes$}}
\put(9.7,1.7){{$\otimes$}}


\put(9,1.81){\line(1,0){0.74}}
\put(10,1.81){\line(1,0){0.74}}
\put(11,1.81){\line(1,0){0.74}}

\put(8.1,1.3){$\scriptstyle 1$}
\put(9.1,1.3){$\scriptstyle 1'$}
\put(10.1,1.3){$\scriptstyle 2$}

\put(8.7,1.7){{$\otimes$}}
\put(9.7,1.7){{$\otimes$}}
\put(10.7,1.7){{$\otimes$}}\put(11.87,1.81){\circle{0.27}}%

\put(3.7,0.82){{$\otimes$}}
\put(4.87,0.93){\circle{0.27}}
\put(4.2,2.58){$\otimes$}
\put(2.95,1.7){{$\otimes$}} 
\put(5.46,1.7){$\otimes$}
\put(5.53,1.94){\line(-5.2,3){1.09}}
\put(3.19,1.92){\line(5,3.2){1.06}}
\put(3.75,1.01){\line(-4,5){0.54}}
\put(4.95,1.01){\line(4,5){0.54}}
\put(3.3,1.07){$\scriptstyle 1$}
\put(4.3,.42){$\scriptstyle 2'$}
\put(5.3,1.07){$\scriptstyle 3'$}
\put(3.17,2.28){$\scriptstyle 1'$}
\put(5.27,2.28){$\scriptstyle 2$}
\put(4.0,.93){\line(1,0){0.74}}
\put(9,1.81){\line(1,0){0.74}}
\put(10,1.81){\line(1,0){0.74}}
\put(11,1.81){\line(1,0){0.74}}

\put(8.1,1.3){$\scriptstyle 1$}
\put(9.1,1.3){$\scriptstyle 1'$}
\put(10.1,1.3){$\scriptstyle 2$}

\put(8.7,1.7){{$\otimes$}}
\put(9.7,1.7){{$\otimes$}}

\end{picture}
\Ec 
\vspace{1cm}

\Bc
\setlength{\unitlength}{0.9cm}
\begin{picture}(10,1)(-0.5,1.70)
\thinlines
  \linethickness{.09mm}
\put(-2.47,1.81){\circle{0.27}}
\put(-1.7,1.7){{$\otimes$}}
\put(-.62,1.7){$\otimes$}\put(.36,1.7){$\otimes$}
\put(-2.4,1.81){\line(1,0){0.74}}
\put(-1.4,1.81){\line(1,0){0.74}}
\put(-0.4,1.81){\line(1,0){0.74}}
\put(-3.1,1.3){$\scriptstyle 1'$}
\put(-2.1,1.3){$\scriptstyle 2'$}
\put(-1.1,1.3){$\scriptstyle 1$}
\put(.1,1.3){$\scriptstyle 3'$}
\put(1.1,1.3){$\scriptstyle 2$}
\put(3.7,0.82){{$\otimes$}}
\put(4.7,0.82){{$\otimes$}}
\put(4.2,2.58){$\otimes$}
\put(3.07,1.87){\circle{0.27}} 
\put(5.46,1.64){$\otimes$}
\put(5.53,1.94){\line(-5.2,3){1.09}}
\put(3.19,1.92){\line(5,3.2){1.06}}
\put(3.75,1.01){\line(-4,5){0.54}}
\put(4.95,1.01){\line(4,5){0.54}}
\put(3.3,1.07){$\scriptstyle 2'$}
\put(4.3,.42){$\scriptstyle 1$}
\put(5.3,1.07){$\scriptstyle 3'$}
\put(3.17,2.28){$\scriptstyle 1'$}
\put(5.27,2.28){$\scriptstyle 2$}
\put(4.0,.93){\line(1,0){0.74}}

\put(9,1.81){\line(1,0){0.74}}
\put(10,1.81){\line(1,0){0.74}}
\put(11,1.81){\line(1,0){0.74}}

\put(8.1,1.3){$\scriptstyle 1$}
\put(9.1,1.3){$\scriptstyle 1'$}
\put(10.1,1.3){$\scriptstyle 2'$}
\put(11.27,1.3){$\scriptstyle 2$}
\put(12.13,1.3){$\scriptstyle 3'$}%

\put(8.7,1.7){{$\otimes$}}
\put(9.87,1.81){\circle{0.27}}
\put(10.7,1.7){{$\otimes$}}
\put(11.76,1.7){$\otimes$}

\end{picture}
\Ec 
\vspace{1cm}
}\eexa

\bexa \label{E5}{\rm We give a simple example to illustrate the difference between local and global names.  Let $\fgo=\fgl(3|4)$.   The anti-distinguished set of simple roots corresponding  to the shuffle $\gs=(1', 2',3',4',1,2,3)$
is 
$$(\gd_1 - \gd_ 2,\; \gd_2 - \gd_3 ,\; \gd_3 - \gd_4 ,\; \gd_4 - \gep_1 ,\; \gep_1 - \gep_2 ,\; \gep_2 - \gep_3).$$  If we add the extending root $\overline{\delta}  + \gep_3 - \gd_1$ and remove $\gep_2 - \gep_3.$
we obtain the set of simple roots with global name 
\be \label{ggb} (\overline{\delta}  + \gep_3 - \gd_1,\; \gd_1 - \gd_ 2,\; \gd_2 - \gd_3 ,\; \gd_3 - \gd_4 ,\; \gd_4 - \gep_1 ,\; \gep_1 - \gep_2)\ee for the Borel with local name $\fbo(\overline{\gs},1 )$, compare \eqref{f19}.  (The new extending root is $\gep_2 - \gep_3.$)  Observe however that if we ignore $\overline{\delta} $ and look at the subscripts in \eqref{ggb} we obtain the permutation $w = (3,1,'2',3',4',1,2)$ which is not a shuffle. The corresponding shuffle is $\overline{\gs}= \nu w$ as in the proof of  Lemma \ref{ymn}.
 
}\eexa

\subsection{Reflection complete Verma modules}\label{sscbs}
We recall the significance of odd reflections in representation theory.
Let  $\fk = \fgl(n|n)$ or $ \fsl(n|m)$ or the affinizations of these algebras.
Let $\fb, \fb'$ be adjacent Borel subalgebras of $\fk$ and suppose 
\be \label{fggb} \fk^\gb \subset \fb, \quad \fk^{-\gb} \subset \fb'\ee for some isotropic root $\gb$. The following is well-known.
\bl \label{cobs} Suppose $V = U(\fk)v_{\lambda}$ where $v_{\lambda}$ is a highest
weight vector for $\fb$ with weight $\lambda$, and set  $u = e_{- \beta} v_{\lambda}$. Then either $u = 0$ or
$u$ is a highest weight vector of weight $\lambda - \beta$ for
$\fb'$. Also $u$ generates a proper 
submodule of $V$ iff 
$(\gl-\gb, \gb) =0$. 
\el \noi 
Consider  a  sequence
 \be \label{distm} 
\mathfrak{b}^{(0)}, \mathfrak{b}^{(1)}, \ldots, \mathfrak{b}^{(r)}\ee
  of Borel subalgebras such that $\mathfrak{b}^{(i-1)}$ and $\mathfrak{b}^{(i)}$ are adjacent for $1 \leq i \leq r$ and $\mathfrak{b}^{(r)}=\fbo' $. 
There are positive odd roots $\gb_i$  of $\fbo$ such that
	\be \label{bar}\fk^{\gb_i} \subset \fb^{(i-1)}, \quad \fk^{-\gb_i} \subset \fb^{(i)}\ee
	for $i\in [r]$.  Without loss we may assume   the roots $\gb_1,\ldots,\gb_r$ are distinct.
	Let $M(\gl)$ be the Verma module induced from the one-dimensional representation of the Borel subalgebra $\fbo$ with weight $\gl$ and highest weight vector $v_\gl$. For   $i \in [r]$ let $x_i$ be  a root vector with weight $-\gb_i$ 
	\be \label{tom}v_0 = v_\gl, \;\;\gl_{0} = \gl, \;\; v_i= x_iv_{i-1},\;\; \gl_{i} = \gl_{i-1}-\gb_i
.\ee
 Then  $v_i$ is a highest weight vector for the Borel subalgebra $\fb^{(i)}$ with weight $\gl_{i}$.  
We say that $M(\gl)$ is {\it reflection complete} if for any sequence of Borel subalgebras as in 
\eqref{distm}, $v_r\neq 0$ 
and 
$U(\fk)v_{r-1}$ is a proper submodule of  $U(\fk)v_{r}.$ By repeated application of Lemma \ref{cobs}, we have
\bc \label{c8} The Verma module $M(\gl)$ is reflection complete if and only if, for all 
chains of Borel subalgebras as in \eqref{distm} we have $v_r\neq 0$ and $\sum_{i=1}^r(\gl-\gb_i , \gb_r) = 0$.\ec\noi 
Now suppose that $\fk = \fgl(n|m)$ or $ \fsl(n|m)$ and 
$\mathfrak{b}^{(0)}=\fbo$ is the distinguished Borel subalgebra. Let 
$\ga= \gep_p-\gd_q$
and using the notation \eqref{323}, let $\mathfrak{b}^{(r)}
\in {\stackrel{{\rm o}}{\cB}}(\ga).$

\bt \label{t8} The Vema module $M(-\gr)$ is reflection complete.
\et\noi 
\bpf
We consider a Levi subalgebra $\bar \fg$ of $\fg$ with highest odd root $\ga= \gep_p-\gd_q$. 
The positive even and odd roots $\bar \Delta_0 $ and $\bar \Delta_1 $  of $\bar \fg$ are given respectively by  
\begin{equation}
\bar \Delta_0=\left\{\gep_{i}-\gep_{j}|p \le i < j\le n\right\}\cup
\left\{\delta_{i}-\delta_{j}|1 \le i < j\le q\right\}.
\notag\end{equation}
and 
\begin{equation}  \label{q3} 
\bar \Delta_1=\left\{\gep_{i}-\delta_{j}|i =p,\dots ,n,\; j=1,\dots q\right\}.
\end{equation}
Let $\gr_i(\bar \Delta)= \frac{1}{2}\sum_{\ga\in  \bar \Delta_i}$ for $i=0,1$ and $\gr(\bar \Delta) = \gr_0(\bar \Delta)-\gr_1(\bar \Delta)$.  When $p=1$ and $q=m$, we denote 
$\gr_i(\bar \Delta)$ and $\gr(\bar \Delta)$ simply by 
$\gr_i$ and $\gr$.  Since 
\be \label{y9} 2\gr_1(\bar \Delta) =q\sum_{i=p}^{n}\gep_{i} - (n+1-p)\sum_{i=1}^{q}\delta_{j}\ee we have 
\be \label{y5} (2\gr_1(\bar \Delta),\ga) =p+q-n-1.\ee In particular
\be \label{y55} (2\gr_1,\ga) =m-n.\ee
Also from  $h_\ga = \sum_{i=p}^{n+q-1}h_i$, it follows that
\be \label{y6}
(\gr,\ga) = n+1-p-q. \ee
Consider a chain of Borel subalgebras as in \eqref{distm} with $\gb_r =\ga$. 
 If $\mathfrak{b}^{(r)}=\fbo(\gs)$ for some $ \gs \in \Sh$ and $\ga= \gep_p - \gd_q$, this means that $p$ and $q'$ are consecutive entries in $\gs$.  Hence the set of roots $\{\gb_i\}_{i=1}^r$ contains all roots $\bar \Delta_1$ from  \eqref{q3}, but does not contain any roots of the form $\gep_i -\gd_q$ or 
$\gep_p -\gd_j$ with $i<p$ or $j>q$.  Thus $ \sum_{i=1}^r\gb_i = 2\gr_1(\bar \Delta)+\gb'$ where $(\gb',\ga)=0$, so by  \eqref{y5}
 and \eqref{y55}
\be \label{q8} \sum_{i=1}^r(\gb_i , \ga) = p+q-n-1.\ee
	Now let $M(-\gr)$ be the Verma module induced from the one-dimensional representation of the Borel subalgebra $\fbo$ with weight $-\gr$ and highest weight vector $v_-\gr$. If $x_i$ are root vectors as in \eqref{tom}, then 
the product $x=x_{r}\ldots x_2x_1$ in  $U(\fgo)$ is non-zero by the PBW theorem, since the weights $\sum_{i=1}^{r} \gb_i $ are distinct and hence $0\neq v_r \in M(-\gr)^{-\gr-\sum_{i=1}^{r} \gb_i} $.  From \eqref{y6} and  \eqref{q8} it follows that $(-\gr-\sum_{i=1}^r\gb_i , \gb_r) = 0$.
Thus the result follows from Corollary \ref{c8}.
\epf\noi 
Now let $\fb$ be the Borel subalgebra of $\Lgh$ extended from the distinguished Borel subalgebra $\fbo$ of $\fgo$. We consider Verma modules induced from $\fb$. 
Since $c$ is central in $\Lgh$, $c$ acts on any Verma module 
$\hat M(\gl)$ as a scalar. If $(c-h^\vee)\hat M(\gl)=0$ where
$h^\vee$ is  the dual Coxeter number, 
we say that $\hat M(\gl)$ is  a Verma module at the {\it critical level.}  
If $\fgo=\fsl(n|m)$ or $\fgo=\fgl(n|m)$ with $m>n$, we have $h^\vee = m-n$, \cite{KaWa} Section 1. Suppose that
$\gl \in \fh^*$ satisfies $\gl(c) = m-n$ and the restriction of $\gl$ to $\fho$ coincides with $-\gr$. Let $\hat M(\gl)$ be  the Verma module for $\Lgh$  induced from the one-dimensional representation of the Borel subalgebra $\fb$ with weight $\gl.$
We have 
\bt \label{t81} The Verma module $\hat M(\gl)=0$  
is reflection complete.
\et\noi
\bpf The proof is similar to the proof of Theorem \ref{t8} so we give fewer details.  Let $x_i$ be root vectors as in \eqref{tom}.  
Suppose $ k\overline{\delta} + \ga$ is a simple root of the 
Borel subalgebra ${\stackrel{{\rm o}}{\fb}}{}(\gs,k)$ of 
${\stackrel{{\rm o}}{\fg}}(k)$ with corresponding root vector $t^{k}e_\ga$. 
By  \eqref{y5}  $(2\gr_1,\ga) =m-n$.  Also 
\be \label{q5} 
[t^{k}e_\ga,t^{-k}e_{-\ga}]= [e_\ga,e_{-\ga}]+ k(e_\ga,e_{-\ga})c ) =h_\ga + kc
\ee
 using \eqref{mjr}. We have  $ \sum_{i=1}^r\gb_i = 2\gr_1(\bar \Delta)+2k\gr_1 +\gb'$ where $(\gb',\ga)=0$, so by  \eqref{y5}
\be \label{q81} \sum_{i=1}^r(\gb_i , \ga) = k(m-n)+p+q-n-1.\ee
Therefore $$(h_\ga +kc)v_r =  x(h_\ga +kc-k(m-n) + p+q-n-1)  v_{\gl} = 0.$$
\epf

\section{Concluding Remarks}\label{ior}
\subsection{Graded posets}\label{grp}

The set $X$ has the structure of a graded poset $X=\bcu_{i=0}^{mn} X_r$ where 
$X_r$ is the set of diagrams containing $r$ boxes.  Then $\sum	_{i=0}^{mn} |X_r|q^r$
is a $q$ binomial coefficient \cite{St} Theorem 6.6.  In Subsection \ref{gax} we defined a graded poset structure on the set of equivalence classes 
$[X\ti \Z]$.  It would be interesting to have any information on the cardinality of $ [X\ti \Z]_r.$

\subsection{A conjecture on infinite orbits}\label{cio}

It seems possible that Young diagrams can be used to study arbitrary orbits of $\cW$ on $\Z^{n|m}$.
The idea to regard a zero entry in $\A(\Gl)$ as an outer corner of a possible path in {\bf R}.  Then using the $W$-action to permute the rows and columns of $\A(\Gl)$, we could try to extend this path as much as possible by adding/deleting outer/inner corners.  Adding or deleting a corner corresponds to applying a morphism to $\Gl$ and the path shows the order in which the morphsims are applied.   These ideas lead to the following conjecture. 

\bco \label{t2} Up to a  translation by a suitable vector $\bA$ 
any infinite orbit  for the $\mathfrak{W}$ action on $\ttk^{n|m}$  has a representative $\Gl_0$  of the form 
 \eqref{e7}.
\eco \noi To illustrate the ideas in the first paragraph, we prove a special case of the conjecture.
Suppose $r_1, \ldots, r_{k}$ and $s_1, \ldots, s_{k}$  are positive integers such that 
\begin{equation}
\label{e19} \sum_{i=1}^k r_i = m \quad \mbox{ and } \quad  \sum_{i=1}^k s_i = n.\ee For $p\in [k-1]$, set
\begin{equation}
\label{w3} a_{p} = \sum_{i=1}^p nr_{i}- ms_{i} 
.\ee 
Note that for $i>1$, 
\begin{equation}
\label{w2} 
nr_{i}- ms_{i} =   a_{i} -a_{i-1}.\ee 
It is convenient to write elements of $\Z^{n|m}$
in block form to match the notation \eqref{prf} used for partitions.  
To do this define for $a\in \Z$  and $j\ge 1$ $$[a]^{j} = (a,a+1,\ldots, a + ({j}-1)) \mbox{ and } ^{j}[a] = (a + ({j}-1)),\ldots,a+1,a).$$ We concatenate such terms.  In  Example \ref{e18} they are  separated by semicolons for clarity. By \eqref{e7}, 
\be \label{crx}{^R\Gl_0}= n([0]^{r_1}, \ldots ,[
(r_1+ \ldots r_{i-1})]^{r_i},\ldots ,[
(r_1+ \ldots r_{k-1})]^{r_k}),\ee
\be \label{e39}^L\Gl_0 =
m({^{s_1}}(s_2+ \ldots s_{k}),\ldots,^{s_i}(s_{i+1}+\ldots+{s_k}), \ldots, ^{{s_k}}0).\ee
\bp \label{prq} Suppose that 
\be \label{e59}{^R\Gl}= n([0]^{r_1}, [a_1]^{r_2},\ldots ,[a_{i-1}]^{r_{i}},\ldots ,
[a_{k-1}]^{r_k})\ee  and 
\be \label{e49}^L\Gl =
m({^{s_1}}[a_1], \ldots, ^{s_i}[a_i], \ldots, ^{s_{k-1}}[a_{k-1}], ^{{s_k}}[0]).\ee
Then $\Gl$ is in the $\cW$-orbit of $\Gl_0$.
\ep 
\bpf 
\noi Below we give   a partition $\gl$ such that 
$^R\Gl={^R\Gl_0} + m\gl'-{\bf mn}$ and $^L\Gl =^L\Gl_0 + n\lst-{\bf mn}$ have the same entries as $^R\Gl$ and 
$^L\Gl$ repectively.  They must then be a permutation of these entries and the result will follow from Corollary \ref{aed}.
In \eqref{vp} and \eqref{vq}  we define index sets $I_i, J_i$ such that there are disjoint unions $[m] = \bcu_{i=1}^k I_i$ and $[n] = \bcu_{i=1}^k J_i$, $|I_i|=r_i$, $|J_1|=s_k$ and $|J_i|=s_{i-1}$ for $i>1$.   We identify  $I_{i}\ti J_{i}$ 
with a sub-rectangle of 
{\bf R} such that $\gep_a-\gd_b \in I_{i}\ti J_{i}$ iff  $b\in I_i$ and $a\in J_i$.  The rectangles are arranged so that $I_1\ti J_1$ fits into the lower left corner of  {\bf R}, and for $t\ge 1$ the rectangle $I_{t+1}\ti J_{t+1}$ lies below and farther to the right than 
$I_{t}\ti J_{t}$, see the diagram below.
These conditions imply that the sets 
$I_t, J_t$ 
have the following form

\be \label{vp}   I_{t} = \{i_{t-1}, \ldots , i_{t}-1\} \mbox{ for } t\in [k]\ee
and 
\be \label{vq}  J_{1} = \{j_{k-1}, \ldots , j_{k}-1\}, \quad J_{t+1} = \{j_{t-1}, \ldots , j_{t}-1\} \mbox{ for }  t\in [k-1]
\ee
where  
$ i_0=j_0=1$, $ i_k=m+1, j_k=n+1$,  
$i_t = r_t +i_{t-1}$ and 
$j_t = s_t +j_{t-1}$.  
Next we draw the path in 
{\bf R}  supporting the rectangles $I_t\ti J_t$ for $t\ge 2$.  The diagram below illustrates the case $k=3$.  The path is drawn with bold lines.  The rectangle is augmented by the vectors $^L\Gl$ and $^R\Gl$ but we show only the lowest (resp.  leftmost) entry in each block of  $^L\Gl$ (resp. $^R\Gl$). 
\Bc
\setlength{\unitlength}{1.2cm}
\begin{picture}(10,4)(-1,-1.00)
\thinlines
  \linethickness{0.05mm}
\put(3.15,2.45){$\scriptstyle 0$}
\put(4.15,2.45){$\scriptstyle a_1$}
\put(5.25,2.45){$\scriptstyle a_2$}
\put(5.25,.35){$\scriptstyle  I_3\ti J_3$}
\put(3.25,-.60){$\scriptstyle  I_1\ti J_1$}
\put(4.25,1.40){$\scriptstyle  I_2\ti J_2$}
\put(2.45,1.25){$\scriptstyle a_1$}
\put(2.45,0.15){$\scriptstyle a_2$}
\put(2.45,-.85){$\scriptstyle 0$}
\multiput(1.9,-1)(1,0){2}{\line(0,1){4.1}}
\multiput(3,-1)(1,0){4}{\line(0,1){4.1}}
 \multiput(1.9,-1)(0,1){4}{\line(1,0){4.1}}
 \multiput(1.9,2.1)(0,1){2}{\line(1,0){4.1}}

  \linethickness{0.5mm}\put(5.0,0){\line(1,0){1}}
\put(4.0,1){\line(1,0){1}}
\put(3.0,2){\line(1,0){1}}
\put(4.0,2){\line(0,-1){1}}
\put(5.0,1){\line(0,-1){1}}
\put(6.0,0){\line(0,-1){1}}
\end{picture}
\Ec 
The partition $\gl$ corresponding  to the path and its dual are given by
\be \label{vd1}  \lst^\trans=((r_1)^{s_{1}},\ldots,(r_1+ \ldots +r_{i})^{s_{i}},\ldots, (r_1+ \ldots +r_{k})^{s_{k}})\ee and
\be \label{vd}  
\gl'=
((s_1+ \ldots +s_{k})^{r_1},\ldots,(s_{i}+\ldots+{s_k})^{r_{i}}, \ldots (s_{k})^{r_k}).\ee
Using \eqref{w3}, the first entry in the $i+1^{st}$ block of the sum $
{^R\Gl_0} + m\gl'$ 
is \by \label{g9} n(r_1+ \ldots +r_{i}) + m(s_{i+1}+
\ldots +s_{k}) 
&=&n(r_1+ \ldots +r_{i}) + m(n -(s_1+ \ldots +s_{i})) \\&=&
 a_i +mn.\ey
Thus the $i+1^{st}$ block of $
{^R\Gl_0} + m\gl' $ is $[a_{i}+mn]^{r_{i+1}}$
The same calculation shows that the $i^{th}$ block of $
{^L\Gl_0} + n\lst$ is   $^{s_{i}}[a_{i}+mn]$.  Since 
$^R\Gl$ and $^L\Gl$ have the same entries  reduced by $mn$, the result follows from Corollary \ref{aed}.

\epf  
\bexa \label{e18}
{\rm \;Let $(n,m)= (7,10), \;(r_1,r_2,r_3)=(6,3,1)$ and $(s_1,s_2,s_3) \;= (3,2,2)$.  Then $a_1= 12, a_2=13$,  
$${^R\Gl_0}= (0,7,14,21,28,35;42,49,56; 63),\quad ^L\Gl_0 =
(60,50,40;30,20;10,0),$$
$${^R\Gl}= (0,7,14,21,28,35;12,19,26;13),\quad ^L\Gl =
(32, 22, 12;23,13;10,0)$$
and $$\gl'=
(7^6,4^3,2^1),\quad \lst^\trans=(10^2, 9^2,6^{3}).$$
} 
\eexa
\br{\rm The significance of the rectangles $I_{i}\ti J_{i}$ is as follows. 
  The zero entries of $\A(\Gl)$ are located in the lower left corners of the rectangles.  Starting with these entries we can use Young diagrams contained in these rectangles in the same way as in Subsection \ref{coL}.}\er
  \bc \label{iol}
The $\cW$-orbit of $\Gl_0$ contains 
$\Gl_1=(0,1,\ldots,n-1|0,1,\ldots,m-1)$.
\ec \noi\bpf
If  $i \in [n]$ let $s_i = 1$, 
 \be \label{izl}R_{i} = \{j \in \Z  |(i-1)m< nj\le im\}\nn\ee
and 
  $r_{i} = |R_{i}|.$  
We have  $|\bcu_{i=1}^pR_{i}| =\lfloor pm/n\rfloor.$  
Hence by \eqref{w3}
\begin{equation}
\label{u3} a_{p} = -pm+\sum_{i=1}^p nr_{i}\in [n]
.\ee 
Also if $a_p \equiv a_q \mod n$, then $p=q$. Hence $\{a_1 ,\ldots, a_n \} =[n]$.  Since $n,m$ are coprime, the set $\{nj|0\le j<m\}=\bcu_{i=1}^n nR_{i}$ is a complete set of residues mod $m$.  We have $nR_{i}-im\subset [m]$, so it follows that $\bcu_{i=1}^n (nR_{i}-im) = [m].$  
\epf

\end{document}